\documentclass[11pt,reqno]{amsart}
\usepackage[text={155mm,225mm},centering]{geometry}
\usepackage[mathscr]{eucal}
\usepackage[all]{xy}
\usepackage{mathrsfs}
\usepackage{amsfonts}
\usepackage{amsmath}
\usepackage{amsthm}
\usepackage{amssymb}
\usepackage{latexsym}
\usepackage{textcomp}
\usepackage{stmaryrd}
\usepackage{hyperref}
\usepackage{color}
\usepackage{graphicx}
\usepackage{float} 
\usepackage{subfigure}
\usepackage{pdfpages}
\newtheorem{theorem}{Theorem}[section]
\newtheorem{lemma}[theorem]{Lemma}

\newtheorem{corollary}[theorem]{Corollary}
\newtheorem{proposition}[theorem]{Proposition}
\newtheorem{proposition-definition}[theorem]{Proposition-Definition}
\newtheorem{example-proposition}[theorem]{Example-Proposition}

\theoremstyle{definition}
\newtheorem{example}[theorem]{Example}
\newtheorem{definition}[theorem]{Definition}
\newtheorem{notation}[theorem]{Notation}

\newtheorem{remark}[theorem]{Remark}

\newtheorem*{convention}{Convention}
\newtheorem*{ack}{Acknowledgements}

\newcommand{\bbZ}{\mathbb{Z}}
\newcommand{\gl}{\mathfrak{gl}}

\newcommand{\mO}{\mathcal{O}}

\newcommand{\GL}{\mathrm{GL}}
\newcommand{\fh}{\mathfrak{h}}

\newcommand{\tos}{\mapsto}

\newcommand{\rmG}{\mathrm{G}}

\newcommand{\ka}{\text{K\"{a}hler}}
\newcommand{\Ker}{ \mathrm{Ker}\,}
\newcommand{\bm}{ \mathbf{m}}
\newcommand{\RA}{_R A}
\newcommand{\bbK}{\mathbb{K}}

\newcommand{\Hom}{\mathrm{Hom}}
\newcommand{\End}{\mathrm{End}}
\newcommand{\DerR}{\mathrm{Der}_{R}}
\newcommand{\Alg}{ {\rm Alg}_R}
\newcommand{\Rep}{\mathrm{Rep}}
\newcommand{\Tr}{\mathrm{Tr} \,}
\newcommand{\bd}{\mathbf{d}}

\newcommand{\Th}{\Theta}

\newcommand{\mD}{\mathcal{D}}

\newcommand{\QQ}{\mathbb{K} \overline{Q}}

\newcommand{\ldb}{\mathopen{\{\!\!\{}}
\newcommand{\rdb}{\mathclose{\}\!\!\}}}
\newcommand{\bw}{\mathbf{w}}

\newcommand{\DDer}{\mathbb{D}\mathrm{er}}

\newcommand{\mG}{\mathcal{G}}

\newcommand{\qneck}{{\mathbb{K} \overline{Q}}_{\hbar}}

\newcommand{\al}{\alpha}

\newcommand{\Spec}{\mathrm{Spec}\,}

\newcommand{\HH}{\mathsf{HH}}
\newcommand{\br}{\mathbf{r}}
\newcommand{\repqd}{\mathrm{Rep}^{Q}_{\mathbf{d}}}
\newcommand{\repqqd}{\mathrm{Rep}^{\overline{Q}} _{\mathbf{d}}}

\newcommand{\mE}{\mathcal{E}}

\newcommand{\ba}{\mathbf{a}}

\newcommand{\bbN}{\mathbb{N}}
\newcommand{\mF}{\mathcal{F}}

\newcommand{\rmH}{\mathrm{H}}
\newcommand{\mT}{\mathcal{T}}

\newcommand{\SN}{\mathrm{SN}}

\newcommand{\tr}{\mathrm{tr}}
\makeatletter

\newcommand{\Rmnum}[1]{\expandafter\@slowromancap\romannumeral #1@}

\allowdisplaybreaks
\title{NONCOMMUTATIVE HAMILTONIAN STRUCTURES AND QUANTIZATIONS ON PREPROJECTIVE ALGEBRAS}

\author{Hu Zhao}

\address{Institute For Advanced Study In Mathematics, Zhejiang University, Hangzhou, P.R. China}

\email{zhaohumath@zju.edu.cn}

\date{\today}

\begin{document}

\begin{abstract}
	Given a noncommutative Hamiltonian space $A$,
	we prove that the  conjecture ``{\it quantization commutes with reduction}''  holds for $A$.
	We further construct a semidirect product algebra $A \rtimes \mG^A$,
	and establish a correspondence between equivariant sheaves on the representation space and left $A\rtimes\mG^A$-modules.
	In the quiver setting,
	using the quantum and classical trace maps,
	we establish the explicit correspondence between quantizations of a preprojective algebra 
	and those of a quiver variety.

	
\end{abstract}

\maketitle

\setcounter{tocdepth}{2}
\tableofcontents

\section{Introduction}\label{sec: intro}

Recent developments have demonstrated that many important spaces are naturally described within the framework of quiver theory (\cite{Soi2019Lec}).
Quiver varieties have therefore attracted growing interest across numerous research areas.
A fundamental challenge is to characterize quantizations of   general quiver varieties and to establish the conditions under which the localization theorem holds.
Significant progress has been made recently, for example, in the work of Bezrukavnikov and Losev \cite{BL2021Eti}.
However, quantizations of  general quiver varieties are still not fully understood, except in special cases such as type $\widetilde{A}$.
For further details, see \cite{EG2002,EGGO2007HC,BPW2016Qua,BLPB2016II,Los2021Loc}.

On the other hand,
the celebrated Kontsevich–Rosenberg principle (\cite{KonRos2000}) states that a geometric structure on a noncommutative algebra is meaningful precisely when it induces the corresponding classical structure on its representation spaces.
Let $\bbK$ be an algebraically closed field of characteristic~zero.
Let $Q$ be a finite quiver.
We double $Q$ by adding an reverse arrow for each arrow, yielding the quiver $\overline{Q}$.
As shown by Kontsevich \cite{Kontsevich} (see also \cite{CBEG2007,Gin2001}), the path algebra $\QQ$ serves as the ``noncommutative cotangent space'' of $\bbK Q$.
$\QQ$ is endowed with a noncommutative Poisson structure (Proposition \ref{prop: induced Lie bracket}) together with a noncommutative moment map (Definition \ref{def: nc moment map})  
\(\bw = \sum_{a\in Q} (a\,a^{\ast} - a^{\ast}a).\)
As shown in \cite{CBEG2007, Van2008Double}, the Hamiltonian reduction of $(\QQ,\bw)$ yields the preprojective algebra $\Pi Q$.

On the other hand,
Schedler \cite{Sch2005} constructed a $\bbK[\hbar]$-algebra~$\QQ_{\hbar}$ quantizing the noncommutative Poisson structure on~$\QQ$. 
See Definition \ref{def: noncom quantization} and Remark \ref{rk: nc quant} for the quantization in the noncommutative setting.

Therefore, for a given preprojective algebra, constructing the required quantizations of its noncommutative Poisson structure is equivalent to verifying the conjecture that ``quantization commutes with reduction'' in the noncommutative setting.
In other words,
one needs to complete the diagram
\[
\begin{aligned}
	\xymatrixcolsep{4pc}\xymatrix{
		\QQ_{\hbar}\ar@{-->}[r]^-{\rm reduction}& ?\\
		\QQ\ar@{~>}[u]^-{\rm quantization}\ar@{-->}[r]_-{\rm reduction}&\Pi Q\ar@{~>}[u]_-{\rm quantization}
	}
\end{aligned}
\]
and clarify its analogue for quiver representation spaces:
\[
\begin{aligned}
	\xymatrixcolsep{4pc}\xymatrix{
		\mathcal{D}_{\hbar}(\Rep^Q_\bd)\ar@{-->}[r]^-{\rm reduction}&
		\displaystyle\frac{\bigl(\mathcal{D}_{\hbar}(\Rep^Q_\bd)\bigr)^{\gl_\bd(\bbK)}}%
		{\bigl(\mathcal{D}_{\hbar}(\Rep^Q_\bd)(\tau - \hbar\chi)(\gl_\bd(\bbK))\bigr)^{\gl_\bd(\bbK)}}\\
		\bbK[T^{*}\Rep^Q_\bd]\ar@{~>}[u]^-{\rm quantization}\ar@{-->}[r]_-{\rm reduction}&\bbK[\mu^{-1}(0)//\GL_\bd(\bbK)]. \ar@{~>}[u]_-{\rm quantization}
	}
\end{aligned}
\]
See Section~\ref{subsec: [Q,R]=0 on rep} for details.

The main results are summarized as follows.
Let $R=\bigoplus_{i\in I}\bbK e_i$ be the commutative ring generated by pairwise orthogonal idempotents $\{e_i\}$. Suppose $A$ is an  $R$-algebra endowed with a double Poisson bracket $\ldb-,-\rdb$ and a moment map $\bw\in A$.
Let $\mG^A$ be the noncommutative gauge group (Definition \ref{def: nc gauge group}).
The Hamiltonian structure is reformulated as the two-term complex
\[
\xymatrix{
	0\ar[r]&\mG^A\ar[r]^-{\xi}&A\ar[r]&0,
}
\]
whose cohomology gives the noncommutative Hamiltonian reduction.

Our first main result is the construction of a semidirect product algebra $A\rtimes\mG^A$ and the correspondence between equivariant sheaves on the representation space and left $A\rtimes\mG^A$-modules.

\begin{theorem}[Theorem \ref{thm: equ.vdb}]
	Let $(A,\ldb-,-\rdb,\bw)$ be a noncommutative Hamiltonian space. For any $N\in\bbN$ and $\GL_N(\bbK)$-equivariant $\mO$-module $\mF$, the sheaf
	\(\mE_N\otimes_{\mO}\mF\)
	naturally carries a  left $A\rtimes\mG^A$-module structure.
\end{theorem}

Our second main result is to develop a combinatorial way
to construct a complex $\widehat\Lambda_{A,\br}^\bullet$ (with parameter $\br\in R$) whose cohomology quantizes the noncommutative Hamiltonian reduction.

\begin{theorem}[Theorem \ref{thm: quantization by a complex}]
	Let $(A,\ldb-,-\rdb,\bw)$ be a noncommutative Hamiltonian space.
	Let $A_\hbar$ quantize $A$ in the sense of Definition \ref{def: noncom quantization}.
	Then for any $\br\in R$, $H^1(\widehat\Lambda_{A,\br}^\bullet)$ is a quantization of $H^1(\Lambda_A^\bullet)$.
\end{theorem}

This establishes the conjecture “{\it quantization commutes with reduction}” in general:

\[
\xymatrixcolsep{4pc}\xymatrix{
	A_\hbar\ar@{-->}[r]^-{\rm reduction}&H^1(\widehat\Lambda_{A,\br}^\bullet)\\
	A\ar@{~>}[u]^-{\rm quantization}\ar@{-->}[r]_-{\rm reduction}&H^1(\Lambda_A^\bullet). \ar@{~>}[u]_-{\rm quantization}
}
\]
Since our construction uses Hochschild homology, all results are Morita invariant.

Our third main result is to show that, in the quiver setting, the quantum trace map yields the correspondence between noncommutative quantum moment maps and quantum moment maps on representation spaces.
\begin{theorem}[Theorem \ref{thm: noncom quant moment map}]
	Let $Q$ be a finite quiver. Let $\bd$ be a dimension vector.
	Then the map
	\begin{equation*}
		\gl_{\bd}(\bbK) \to \mD_{\hbar}(\repqd),\ v \tos \tr([\widehat{\bw}] v)
	\end{equation*}
	is a quantum moment map for $\repqqd$.
	Furthermore, for an arbitrary $\br \in R$,
	\begin{equation*}
		\gl_{\bd}(\bbK) \to \mD_{\hbar}(\repqd),\ v \tos \tr\Big(([\widehat{\bw}] + \hbar \sum_{i \in Q_0} r_i I_i ) v \Big)
	\end{equation*}
	is also a quantum moment map,
	$I_i$ is the identity matrix at $i$-th component and zero elsewhere.
\end{theorem}

Another notable property is that there is an explicit correspondence between
noncommutative quantizations on the preprojective algebra
and those on the quiver variety.

\begin{lemma}[Lemma \ref{lem: q.trace preserves ideals}]
	Let $Q$ be a finite quiver.
	Let $\mathbf{d}$ be a dimension vector.
	For an arbitrary $\br \in R$, there is a unique character 
	$\chi_{\br}$ of $\gl_{\bd} (\bbK)$
	such that
	\begin{equation*}
		\mathrm{Tr}^{q}\Big(\qneck \{ \ell_{\hbar}[p \bw]+ \ell_{\hbar}([p\br])\hbar \, \vert \, [p \bw] \in (\QQ \bw \QQ)_\natural \}\qneck \Big)\end{equation*}
	is contained in
	\(\displaystyle{\Big( \mathcal{D}_{\hbar}(\Rep^Q _\bd) (\tau - \hbar \chi_{\br}}) 
	(\gl_{\bd} (\bbK)) \Big)^{\gl_{\bd} (\bbK)}.\)
\end{lemma}

The character $\chi_{\br}$ is given in an explicit way:
\(\displaystyle{\chi_{\br}=\sum_{k \in Q_{0}}\Big(- \sum_{a\in Q,s(a)=k} d_{t(a)} +r_k\Big)tr_{k}.}\)

Then, it is clear that the noncommutative version of ``{\it quantization commutes with reduction}'' fits into the Kontsevich–Rosenberg principle
via the commutative cubic (Proposition \ref{prop: [tr,quant]}, Theorem \ref{thm: quant cond after red})
\begin{equation*}\label{comcubic}
	\begin{split}
		\xymatrixrowsep{0.8pc}
		\xymatrixcolsep{1.2pc}
		\xymatrix{
			\qneck \ar@{-->}[rd] \ar[rr]^-{\mathrm{Tr}^{q}}  && 
			\mD_{\hbar}(\Rep^Q _{\bd}) \ar@{-->}[rd] \\
			&\Pi Q_{\hbar, \br} \ar[rr]^{\mathrm{Tr}^{q}}  
			&& 
			\mathcal{M}_{\bd}(Q)_{\hbar, \chi_{\br}} \\
			\mathbb{K}\overline{Q} \ar@{~>}[uu] \ar@{-->}[rd] \ar'[r][rr]^-{\mathrm{Tr}}&& 
			\mathbb{K}[T^\ast \Rep^Q _\bd] \ar@{~>}'[u][uu] \ar@{-->}[rd] \\
			& \Pi Q \ar@{~>}[uu] \ar[rr]^{\mathrm{Tr}} 
			&& \mathbb{K}[\mathcal{M} _{\mathbf{d}}(Q)]. \ar@{~>}[uu].
		}
	\end{split}
\end{equation*}

At the end,
we highlight the following novel contributions compared to \cite{Zhao2021Com}.
\begin{enumerate}
	\item In  \cite{Zhao2021Com},
	the quantum preprojective algebra is constructed only for $\br = 0$.
	In this work, our complex formalism naturally explains the role of the parameter $\br$.
	Only when the noncommutative Hamiltonian reduction is realized as a complex and the correspondence between noncommutative fields and Hamiltonians is remembered by the map $\xi$,
	then correspondence between quantum fields and quantum Hamiltonians canonically contains higher-order information.
	In this case,
	those at order $1$ is decoded by $\br$.
	
	\item This work establishes a constructive framework applicable to deformed preprojective algebras (Section \ref{subsec: deformed case}).
\end{enumerate}

The structure of this paper is as follows.


In Section \ref{sec: noncom Ham sys},
we recall noncommutative Hamiltonian geometry.
Then, 
the semidirect product algebra $A \rtimes \mG^A$ is introduced and  its representations are related to  equivariant sheaves on the representation space.
Finally, the conjecture ``{\it quantization commutes with reduction}'' in the noncommutative setting is proved.

In Section \ref{sec: quiver}, we recall the noncommutative Hamiltonian structure on the path algebra $\QQ$ associated with a quiver $Q$. We show that the quantum trace maps are preserved under quantum reduction. In particular, we establish a correspondence between quantizations of preprojective algebras and those of quiver varieties. Finally, we show that the noncommutative ``{\it quantization commutes with reduction}'' fits into the Kontsevich--Rosenberg principle.

\begin{ack}
	The author is deeply grateful to Professor Xiaojun Chen and Professor Farkhod Eshmatov for their insightful discussions.
	Special thanks are extended to Professor Yongbin Ruan for his continuous support and valuable advice.
	In particular, the author would like to thank Professor Hiraku Nakajima who suggests the quantization-localization problem for noncommutative quantization.
	The author has greatly benefited from discussions with Professor Christopher Brav and Professor Si Li.
	This research was funded by the Postdoctoral Fellowship Program of CPSF(GZC20232337).
\end{ack}

\begin{convention}
	\begin{itemize}
		\item $\bbK$ is an algebraically closed field of characteristic zero.
		
		\item Throughout this work, an algebra refers to a finitely generated algebra, not necessarily commutative.
		
		\item Throughout this work, a module refers to a finitely generated module.
		
		
		
		\item For a quiver $Q$, the vertex set is denoted by $Q_0$, and the set of arrows is still denoted by $Q$.
		For an arrow $a$, $s(a)$ is the source of $a$, $t(a)$ is the target of $a$.
		
		
	\end{itemize}
\end{convention}

\section{Noncommutative Hamiltonian spaces}\label{sec: noncom Ham sys}

In this section, we recall counterparts of Hamiltonian spaces in the noncommutative setting.

\subsection{Noncommutative Hamiltonian spaces}\label{subsec: double Poi}

Throughout this work, we consider algebras over a commutative ring $R$ constructed as follows.
Let $I$ be a finite index set. The commutative ring $R$ is defined as 
$\bigoplus_{i \in I} \bbK e_i$, where $e_i$ are pairwise orthogonal idempotents.
We abbreviate $-\otimes_R -$ to $- \otimes -$.

Firstly, we  recall representation spaces.
Let \(V\) be a finite-dimensional \(R\)-module over \(\bbK\). 
The {\it representation space of an \(R\)-algebra \(A\) on \(V\)} is
defined as 
\(\Rep^A_V :=\;\Hom_{\Alg}\bigl(A,\gl(V)\bigr),\)
parametrizing left \(A\)-module structures on \(V\).
Moreover, the general linear group \(\GL(V)\) acts on \(\Rep^A_V\) by conjugation:
for \(\rho\in\Rep^A_V\), \(g\in\GL(V)\), and \(a\in A\),
\((g.\rho)(a) : =g\,\rho(a)\,g^{-1}.\)

In this work,
we restrict to $V=\bbK^\bd$
where $\bd = (d_i)$ is a dimension vector in $\bbZ^{I} = \oplus_{i \in I} \bbZ$.
The associated representation space is denoted by $\Rep^A_\bd$.
Each \(a\in A\) defines a \(\gl(\bbK^\bd)\)-valued function
\begin{equation}\label{matrix function}
	(a)\colon \Rep^A_\bd \;\to\;\gl(\bbK^\bd),\quad
	\rho\mapsto \rho(a).
\end{equation}
If \(a\in e_jAe_i\), then \(\rho(a)\) corresponds to  a \(d_j\times d_i\) matrix, whose \((p,q)\)-entry is denoted by \(\rho(a)_{p,q}\).
This yields a regular function
\[
(a)_{p,q}\colon \Rep^A_\bd\;\to\;\bbK,\quad
\rho\mapsto \rho(a)_{p,q},
\]
which is often abbreviated to \(a_{p,q}\).
The coordinate ring $A_\bd : = \bbK[\Rep^A_\bd]$ is generated by the matrix coefficients
$$
\bigl\{ a_{p,q} \bigm| a \in A,\ 1 \leq p \leq d_{t(a)},\ 1 \leq q \leq d_{s(a)} \bigr\}.
$$
See \cite[Proposition.~12.1.6]{Gin2005non} for a proof.

In the noncommutative setting,
the noncommutative $\ka$ differential is defined as follows.

\begin{definition}\label{def: ka differential}
	Let $A$ be an $R$-algebra.
	The {\it noncommutative K\"ahler differential of $A$}  is defined to
	be the $A$-bimodule $\Omega^1 \RA :=\Ker (\bm)$.
\end{definition}
Here, $\bm$ denotes the multiplication  $A \otimes  A \to A$.
Furthermore, the noncommutative $\ka$ differential is related to the derivation space via the following standard fact.
\begin{proposition}\label{prop: derhom}
	Suppose $A$ is an $R$-algebra and $M$ is a left $A^e$-module.
	Then there is a canonical bijection
	\(\DerR(A, M) \cong \Hom_{A^e} (\Omega^1 \RA, M).\)
\end{proposition}
Here,
$A^e = A \otimes A^{op}$;
we implicitly use the equivalences among $A$-bimodules, left $A^e$-modules and right $A^e$-modules.

Next, we recall noncommutative vector fields.
Following \cite{CBEG2007}, the ``noncommutative vector fields on $A$'' are defined to be double derivations on $A$.
\begin{definition}
	Let $A$ be an $R$-algebra.
	A {\it double derivation on $A$ }is defined as
	an $R$-derivation from  $A$ to $A \otimes A$,
	where $A \otimes A$ is endowed with the outer $A$-bimodule structure.
\end{definition}
Recall that there are two $A$-bimodule structures on $A\otimes A$:
for any $a,\ b,\ x,\ y$ in $A$, 
\begin{equation}\label{for: outer bimod}
	a(x\otimes y)b : = ( ax)\otimes (yb)
\end{equation}
gives the {\it outer bimodule structure};
\begin{equation}\label{for: inner bimod}
	a \bullet (x\otimes y) \bullet b := (xb)\otimes (ay)
\end{equation}
gives the {\it inner bimodule structure}.
The set of double derivations on $A$ is denoted as
$\mathbb{D}\mathrm{er}_{R}A$.
The inner $A$-bimodule structure endows
$\mathbb{D}\mathrm{er}_{R}A$ with an $A$-bimodule structure.
In subsequent discussions,
we will use the Sweedler notation:
for a double derivation $\Theta: A\rightarrow A\otimes A$,
\(\Theta(a) = \Theta'(a) \otimes \Theta''(a).\)


\begin{example}
	Consider the free $R$-algebra
	$A=R\langle x,y  \rangle$.
	Then according to (\ref{matrix function}), $x, y$  induce matrix-valued functions:
	\((x)=(x_{i, j}),\ (y)=(y_{i, j}) \in \gl_2(\bbK) \otimes \bbK[\Rep^A_2].\)
	Here, $x_{i,j}$ and $y_{i,j}$ are matrix coefficient functions.
	$dx, dy$ induce matrix-valued $1$-forms on representation space:
	\((dx) = (dx_{i, j}),\ (dy) = (dy_{i, j}) \in \gl_2(\bbK) \otimes \Omega^{1} (\Rep^A _2).\)
	Consequently,
	the noncommutative $2$-form $dx dy$ induces matrix-valued $2$-form:
	\((dx_{i, j})\wedge(dy_{i, j}) = \big(\sum_{k=1} ^2 dx_{i, k}\wedge dy_{k, j}\big).\)
\end{example} 

\begin{example}\label{eg: marix der}
	Let $A$ be a $\bbK$-algebra.
	Let $N$ be a positive integer.
	For any double derivation $\Theta \in \DDer_{\bbK} A$ and any $a \in A$,
	the action of the matrix-valued derivation $(\Theta_{i, j})$ on $A_N$,  where each $\Theta_{i, j}$ is a derivation on $A_N$, is given as follows.
	\begin{equation*}
		\Theta_{i, j}(a_{u, v}) = (\Theta'(a))_{u, j} (\Theta''(a))_{i, v}.
	\end{equation*}
\end{example}

Then, we recall double Poisson brackets and noncommutative Poisson structures.
For a positive integer $n$ and a vector space $V$, the symmetric group $\mathcal{S}_n$ acts on $V^{\otimes n}$: for any $\sigma\in\mathcal{S}_n$ and $v_1\otimes\cdots\otimes v_n\in V^{\otimes n}$,
\(\sigma(v_1\otimes\cdots\otimes v_n)
 : = v_{\sigma^{-1}(1)}\otimes v_{\sigma^{-1}(2)}\otimes\cdots\otimes v_{\sigma^{-1}(n)}.\)


\begin{definition}[Van den Bergh]\label{def: 2-bracket}
	Let $A$ be an $R$-algebra. A \emph{double Poisson bracket} on $A$ is an $R$-bilinear map
	\[
	\ldb-,-\rdb\colon A\otimes A\to A\otimes A,
	\]
	satisfying the following axioms:
	for any $a,b,c\in A$,
	\begin{enumerate}
		\item $\ldb a,b\rdb=-\ldb b,a\rdb^\circ$;
		\item $\ldb a,bc\rdb=\ldb a,b\rdb\,c + b\,\ldb a,c\rdb$;
		\item $\ldb a,\ldb b,c\rdb\rdb_L
		+ (132)\,\ldb c,\ldb a,b\rdb\rdb_L 
		+ (123)\,\ldb b,\ldb c,a\rdb\rdb_L = 0$.
	\end{enumerate}
\end{definition}
Here $\ldb b,a\rdb^\circ=\ldb b,a\rdb''\otimes\ldb b,a\rdb'$;
for $a\in A$ and $b=b_1\otimes\cdots\otimes b_n\in A^{\otimes n}$,
\[
\ldb a,b\rdb_L
=\ldb a,b_1\rdb\otimes b_2\otimes\cdots\otimes b_n.
\]

Van den Bergh (\cite{Van2008Double}) showed that such a double Poisson bracket naturally induces a Lie algebra structure on the zeroth Hochschild homology.
Setting
\[
\{a,b\}: = \bm (\ldb a, b \rdb) = \ldb a,b\rdb'\,\ldb a,b\rdb'' ,
\]
one obtains:

\begin{proposition}\cite[Corollary 2.4.6]{Van2008Double}\label{prop: induced Lie bracket}
	Let $(A,\ldb-,-\rdb)$ be a double Poisson algebra. Then $\HH_0(A)=A/[A,A]$ inherits a Lie algebra structure via the bracket $\{-,-\}$.
\end{proposition}

Throughout this work, that induced Lie bracket $\{-,-\}$ is our notion of \emph{noncommutative Poisson structure}, as it lives on Hochschild homology and is Morita invariant.

A noncommutative analogue of Hamiltonian $\rmG$-spaces in differential geometry
is the notion of a noncommutative Hamiltonian space:
\begin{definition}[Crawley--Boevey--Etingof--Ginzburg, Van den Bergh]\label{def: nc moment map}
	Let $(A,\ldb-,-\rdb)$ be a double Poisson algebra. A \emph{noncommutative moment map} is an element
	\(\bw=\sum_{i \in I}\bw_i\;\in\;\bigoplus_{i\in I}e_iAe_i\)
	such that, for every $p\in A$,
	\[
	\ldb\bw,p\rdb
	=\sum_{i \in I}\bigl(p\,e_i\otimes e_i-e_i\otimes e_i\,p\bigr).
	\]
\end{definition}

\begin{definition}[Crawley--Boevey--Etingof--Ginzburg, Van den Bergh]\label{def: noncom Hamil space}
	A \emph{noncommutative Hamiltonian space} is a double Poisson algebra $(A,\ldb-,-\rdb)$ endowed with a noncommutative moment map $\bw$.
\end{definition}

The compatibility of this construction with the Kontsevich--Rosenberg principle
is guaranteed by the following proposition.
One can find proof in \cite[Proposition 7.11.1]{Van2008Double} and \cite[Theorem 6.4.3]{CBEG2007}.

\begin{proposition}\label{prop: noncom ham to ham}
	Let $(A,\ldb-,-\rdb,\bw)$ be a noncommutative Hamiltonian space. Then for an arbitrary dimension vector $\bd$, 
	$\Rep^A_\bd$ is a Poisson $\GL_\bd(\bbK)$-space, with the moment map
	\[
	\mu\colon\Rep^A_\bd\to\gl_\bd(\bbK)^*,\quad\rho\mapsto\tr\bigl(\rho(\bw)\,\cdot -\bigr).
	\]
	Moreover, the corresponding Poisson bracket on the coordinate ring $A_\bd$ is given by the following formula:
	\begin{equation}\label{for: double poisson to poisson}
		\{a_{ij},b_{uv}\}
		=\ldb a,b\rdb'_{u,j}\,\ldb a,b\rdb''_{i,v},
	\end{equation}
	for $a,b\in A$.
\end{proposition}

As established by Van den Bergh,
the group action on the representation space are related to \emph{gauge elements},
which are double derivations defined by
\begin{equation}\label{for: gauge elements}
	E_i(a) := \ldb \bw_i,a\rdb,
	\quad
	E := \sum_{i \in I} E_i.
\end{equation}

This correspondence is made explicit in the following proposition.
\begin{proposition}\cite[Proposition 7.9.1]{Van2008Double}\label{prop: gauge elements}
	Let $e_{p,q} ^i \in \gl_{\bd} (\bbK)$ denote the  elementary matrix at the $i^{\rm th}$ component,
	which is
	$1$ in the $(p,q)$-entry and zero everywhere else. 
	Then $(E_i)_{p,q}$ (Example \ref{eg: marix der}) acts as $e_{q,p}^i$ on $A_\bd$.
\end{proposition}
\begin{proof}
	Direct calculation
	shows that
	for any $a \in A$,
	\begin{equation*}\label{Eipq}
		(E_i)_{p,q} (a_{u,v})=\delta_{s(a), i}  \delta_{p,v}a_{u,q} -\delta_{t(a), i} \delta_{u,q} a_{p,v} = e^i _{q,p}.a_{u,v}.
	\end{equation*}
\end{proof}

Following the Kontsevich--Rosenberg principle, the noncommutative gauge group is defined as follows.
\begin{definition}\label{def: nc gauge group}
	Let $A$ be an $R$-algebra.
	The {\it noncommutative gauge group $\mG^A$ of $A$} is defined to be the $A$-bimodule generated by gauge elements
	\(\{ E_{i} \vert  i \in I\}.\)
\end{definition}

Finally, we recall noncommutative Hamiltonian reduction.
\begin{definition}\cite{CBEG2007, Van2008Double}
	Let $(A, \ldb-,-,\rdb, \bw)$ be a noncommutative  Hamiltonian space.
	The {\it noncommutative Hamiltonian reduction of $(A, \ldb-,-\rdb, \bw)$}
	is
	the quotient algebra
	\(\displaystyle{\frac{A}{A \bw A}}.\)
\end{definition}

Consider a complex of $A$-bimodules
\begin{equation}\label{for: 2 complex}
	\xymatrix{
		\Lambda_{A} ^{\bullet}: 0 \ar[r]& \mG^A \ar[r]^-{\xi}& A \ar[r]&0.
	}
\end{equation}
Here, $\mG^A$ is of degree $0$;
$A$ is of degree $1$.
$\xi$  assigns $p \bw_i q $ to $p E_{i}q \in \mG^A$.
It is clear that for a noncommutative Hamiltonian space $(A, \ldb-,-\rdb,  \bw)$,
\(\displaystyle{H^0 (\Lambda_{A} ^{\bullet}) = 0\text{ and }\ H^1(\Lambda_{A} ^{\bullet}) = \frac{A}{A \bw A}}.\)

The compatibility of this construction with the Kontsevich--Rosenberg principle
is guaranteed by the following proposition.
\begin{proposition}\label{prop: noncom red}\cite{CBEG2007, Van2008Double}
	Let $(A,\ldb-,-\rdb,\bw)$ be a noncommutative Hamiltonian space. Then:
	\begin{enumerate}
		\item The Lie bracket on $\HH_0(A)$ descends to $\HH_0\bigl(A/A\,\bw\,A\bigr)$, and the projection
		\[
		\HH_0(A)\;\to\;\HH_0\bigl(A/A\,\bw\,A\bigr)
		\]
		is a Lie homomorphism.
		\item The categorical quotient $\Rep^{A/A\,\bw\,A}_\bd//\GL_\bd(\bbK)$ is the Hamiltonian reduction of $\Rep^A_\bd$.
	\end{enumerate}
\end{proposition}

When an algebra $B$ arises from $A$ by (noncommutative) Hamiltonian reduction, we write
\[
\xymatrix{A\ar@{-->}[r]&B.}
\]

\subsection{Equivariant sheaves and the semidirect product algebra}\label{subsec: equi sheaves}

Firstly,
Van den Bergh (\cite[Section 3]{Van2008Double}) introduces a graded double Poisson bracket $\ldb-,-\rdb_{\rm SN}$
on the tensor algebra $T_A \DDer(A)$, whose elements in $\DDer(A)$ are of degree 1.
This double bracket is known as the double Schouten bracket.
Based on his construction,
the semidirect product algebra is defined as follows.
\begin{definition}\label{def: semiproduct alg}
	Let $(A, \ldb-,-\rdb, \bw)$ be a noncommutative Hamiltonian space.
	The {\it  semidirect product algebra $A \rtimes \mG^A$ associated with $(A, \ldb-,-\rdb, \bw)$} is defined as the quotient algebra of the free product $A \cdot \mG^A$ by following relations.
	\begin{enumerate}
		\item For any  $\Th_1,\Th_2 \in \mG^A$,
		$\Th_1 \cdot \Th_2 - \Th_2 \cdot \Th_1 = \{\Th_1, \Th_2\}_{\rm SN} $.
		
		\item For any $a_1, a_2 \in A$,
		$a_1 \cdot a_2 = a_1 a_2$.
		
		\item For any $\Theta \in \mG^A$ and any $a \in A$,
		$\Theta \cdot a - a \cdot \Th = \{\Th, a\}_{\SN}$.
	\end{enumerate}
\end{definition}
From now on, the bracket $\{\Th, a\}_{\SN}$ in $(3)$ will be denoted by $a^\Th$.

In practice, moduli spaces often attract more attention.
A common approach to studying the moduli stack $[\Rep^A_\bd /\GL_{\bd}(\bbK)]$
is via the category of $\GL_\bd(\bbK)$-equivariant sheaves on $\Rep^A_\bd$.
Therefore, we aim to relate $\GL_\bd(\bbK)$-equivariant sheaves on $\Rep^A _\bd$ to  $A \rtimes \mG^A$-modules.

We recall fundamental concepts in equivariant sheaves theory.
Details can be found in \cite{CG2010Rep}.
Let $\rmH$ be a linear algebraic group over $\bbK$.
Let $X$ be an algebraic $\rmH$-variety.
Let $\fh$ be the Lie algebra of $\rmH$.
The action morphism $\ba: \rmH \times X \to X$ induces the infinitesimal action 
$\tau: \fh_X \to \mT_X$, where $\fh_X$ is the localization $\fh \otimes_\bbK \mO_X$, and 
$\mT_X$ is the tangent sheaf of $X$.
An $\rmH$-equivariant structure on an $\mO_X$-module $\mF$ induces a Lie algebra morphism
\(\kappa: \fh \to \End_{\bbK}(\mF)\) such that for $\gamma \in \fh$, $a \in \mO_X$ and $f \in \mF$,
\(\kappa(\gamma)(a f) = \gamma(a)f + a \kappa(\gamma)(f).\)

For simplicity, we consider $V = \bbK^N$.
Let $\mE_N$ be the section sheaf of  the tautological bundle
\(\Rep^A _N \times_{\rm \GL_N(\bbK)} \gl(\bbK^N).\)
For any $\GL_N (\bbK)$-equivariant  $\mO$-module $\mF$ on $\Rep^A_N$, it is clear that
$\mE_N \otimes_{\mO} \mF$
is a sheaf of $A$-bimodules.
Here $\mO$ is the structure sheaf on the representation space.
For  an arbitrary section $(f_{i,j}) \in \mE_N \otimes_{\mO} \mF$
and $\Theta \in \mG^A$,
the action of $\Theta$ on $(f_{i,j})$
is defined by the formula for the entry of $\Theta(f_{i,j})$ at $(u,v)$: 
\begin{equation}\label{for: HC}
	(\sum_{i = 1} ^{N} \Theta_{i,i}). f_{u,v} = (\Tr \Theta). f_{u,v}.%
\end{equation}
Here $\Theta_{i,i}. f_{u,v}$ is given by the equivariant structure   on $\mF$ 
which is precisely the $\kappa$ as above.

\begin{theorem}\label{thm: equ.vdb}
	Let $(A, \ldb-,-\rdb, \bw)$ be a noncommutative Hamiltonian space.
	Then for an arbitrary $N \in \bbN$ and any $\GL_N(\bbK)$-equivariant  $\mO$-module $\mF$,
	\(\mE_N \otimes_{\mO} \mF\)
	naturally carries a   left $A \rtimes \mG^A$-module structure.
\end{theorem}

\begin{proof}
	Since the $A$-action and $\mG^A$-action on $\mE_N \otimes_{\mO} \mF$ is given,
	to prove $\mE_N \otimes_{\mO} \mF$ is a sheaf of left $A \rtimes \mG^A$-module,
	What  needs to be checked is the compatibility with Definition \ref{def: semiproduct alg}.
	
	By (\ref{for: HC}), the entry of $[\Th, \Phi](f_{i,j})$ at $(u,v)$ is actually given by 
	\(\Tr \Th . \Tr \Phi . f_{u,v} - \Tr \Phi . \Tr \Th . f_{u,v} = \{\Tr \Th, \Tr \Phi\}. f_{u,v}.\)
	By \cite[Proposition 7.7.3]{Van2008Double},
	$\Tr$ commutes with Schouten brackets;
	then $\{\Tr \Th, \Tr \Phi\}$ equals to $\Tr \{\Th, \Phi\}_{\SN}$.
	Here, the bracket $\{\Tr \Th, \Tr \Phi\}$ is the Schouten bracket on poly-vector fields on representation spaces.
	
	$(2)$ in Definition \ref{def: semiproduct alg} is canonically compatible with  (\ref{for: HC}).
	Now, compatibility between $(3)$ in Definition \ref{def: semiproduct alg} and (\ref{for: HC}) is as follows.
	By calculation, the $(u,v)$-entry of 
	\((\Theta a) (f_{i,j})\)
	is
	\begin{equation*}
		\begin{split}
			\Big((\Theta a) (f_{i,j})\Big)_{u,v} & = \Tr \Th . \big(a(f_{i,j})\big)_{u,v}\\
			& = \Tr \Th . \big( \sum_{k=1}^{N} a_{u,k} f_{k,v}\big)\\
			& = \sum_{i =1}^N \Th_{i,i} . \big( \sum_{k= 1} ^N a_{u,k} f_{k,v}\big)\\
			& = \sum_{i =1}^N \sum_{k = 1} ^N \Big(\Theta_{i,i}(a_{u,k}) f_{k,v}  + a_{u,k} \Theta_{i,i} (f_{k,v}) \Big)\\
			& = \sum_{i =1}^N \sum_{k = 1} ^N \Big(\Theta'(a)_{u,i}\Th''(a)_{i,k} f_{k,v}  + a_{u,k} \Theta_{i,i} (f_{k,v}) \Big).
		\end{split}
	\end{equation*}
	On the other hand,
	the $(u,v)$-entry of \(	(a^\Th + a\Th) (f_{i,j}) \) is
	\begin{equation*}
		\begin{split}
			\Big((a^\Th + a\Th) (f_{i,j})\Big)_{u,v} & = \sum_{k=1} ^N (a^\Th)_{u,k} f_{k,v} + \big(a. \Th. (f_{i,j})\big)_{u,v}\\
			& = \sum_{k=1} ^N (a^\Th)_{u,k} f_{k,v} + \sum_{k=1} ^N a_{u,k}.\Tr \Th . f_{k,v}\\
			& =   \sum_{k=1} ^N (\{\Th, a \} )_{u,k} f_{k,v} + \sum_{k=1} ^N a_{u,k}.\Tr \Th . f_{k,v}\\
			& = \sum_{i =1}^N \sum_{k = 1} ^N \Big(\Theta'(a)_{u,i}\Th''(a)_{i,k} f_{k,v}  + a_{u,k} \Theta_{i,i} (f_{k,v}) \Big).
		\end{split}
	\end{equation*}
	Thus, $\mE_N \otimes_{\mO} \mF$ is a left $A \rtimes \mG^A$-module.
\end{proof}

\subsection{Noncommutative quantum reduction}\label{subsec: noncom quant red}

Let $(A, \ldb-,- \rdb)$ be a double Poisson algebra.
By
Proposition \ref{prop: induced Lie bracket},
$\HH_{0}(A)$ carries a natural  Lie algebra structure.
Throughout this work, following the idea in \cite{Sch2005} and \cite{GinSch2006},
we regard a quantization of $(A, \ldb-,-\rdb)$ as a PBW-deformation of the induced Lie algebra $\HH_{0}(A)$.
\begin{definition}\label{def: noncom quantization}
	Let $(A, \ldb-,- \rdb)$ be a double Poisson algebra.
	A {\it quantization of the noncommutative Poisson structure on $A$} is a $\bbK[\hbar]$-algebra $A_\hbar$ together with an isomorphism
	$$
	\ell_{\hbar} : Sym(\HH_{0}(A))[\hbar] \to A_{\hbar}
	$$
	of $\bbK[\hbar]$-modules  such that for any $x, y \in \HH_{0}(A)$,
	\begin{equation}\label{noncom dirac pic}
		\ell_\hbar(\{x,y\}) = \frac{-1}{\hbar} [\ell_{\hbar}(x), \ell_{\hbar}(y)]\ {\rm mod}\ \hbar.
	\end{equation}
\end{definition}
In the absence of ambiguity,
$A_{\hbar}$ is called  a quantization of $A$.

\begin{remark}\label{rk: nc quant}
	The Hochschild-Kostant-Rosenberg theorem says that
	Hochschild homology $\HH_{\bullet} (S)$ of a smooth affine scheme $\Spec S$
	is isomorphic to
	the de Rham complex
	$\Omega^{\bullet}_{S}$ and
	$\HH_{0}(S)$ is the algebra of functions.
	Therefore, Definition \ref{def: noncom quantization} is reasonable.
\end{remark}

If $S$ is a quantization of $B$,
we write it as 
\begin{displaymath}
	\xymatrix{
		B  \ar@{~>}[r]              & S.
	}
\end{displaymath}

Recall that  a two-term complex $(\ref{for: 2 complex})$ is constructed for a noncommutative Hamiltonian space $(A, \ldb-,-\rdb,  \bw)$,
such that the  cohomology  gives the noncommutative Hamiltonian reduction.
Analogously, if the noncommutative Hamiltonian space admits a quantization $A_{\hbar}$, it is natural to expect
a new complex whose cohomology  gives a quantization of the noncommutative Hamiltonian reduction.
In other words, the goal is to construct a complex with a parameter $\br \in R$:
\begin{equation*}
	\xymatrix{
		\widehat{\Lambda}_{A,\br} ^{\bullet}: 0 \ar[r] & {\mG^A _{\hbar}} \ar[r] & A_\hbar \ar[r] & 0,
	}
\end{equation*}
such that one has a diagram
\begin{equation*}
	\xymatrix{
		0 \ar[r] & {\mG^A _{\hbar}} \ar[r]  & A_{\hbar} \ar[r] & 0\\
		0 \ar[r]& \mG^A \ar[r]^-{\xi}   \ar@{~>}[u]& A \ar[r]   \ar@{~>}[u]&0 
	}
\end{equation*}
and
\begin{equation}\label{quantization on cohomology}
	\xymatrix{
		H^1({\Lambda}_{A} ^{\bullet})  \ar@{~>}[r]              & H^1(\widehat{\Lambda}_{A,\br} ^{\bullet}).
	}
\end{equation}

Firstly, the quantized noncommutative gauge group $\mG^A _{\hbar}$ is constructed as follows.
The construction is motivated by the correspondence between Hamiltonians and their quantum counterparts.
Let $S$ be the Poisson algebra of smooth functions on a Poisson manifold.
Let $S_\hbar$ be a deformation quantization of $S$.
Fix a Hamiltonian $H \in S$, the Hamiltonian flow is defined by the field $\{H,-\}$.
Then, the corresponding quantum Hamiltonian field, as an endomorphism of the algebra $S_{\hbar}$ of quantum observables, is given by the commutator action ${\rm ad}_{\hat{H}} = [\hat{H},-]$,
where $\hat{H}$ is a lifting of $H$ in $S_\hbar$.

As before, gauge elements $\{E_i = \ldb \bw_i,- \rdb\}$ are the double Hamiltonian vector fields associated with the noncommutative Hamiltonians $\{\bw_{i}\}$.
Therefore,
the noncommutative quantum gauge elements $\{\widehat{E_i}\}$ are constructed via lifting $\{\bw_i\}$ to $A_{\hbar}$.

\begin{definition}
	Let $(A, \ldb-,-\rdb,\bw)$ be a noncommutative Hamiltonian space.
	Let $A_{\hbar}$ be a quantization of A.
	Then {\it quantum gauge elements} are inner derivations
	\[
	\big\{ \frac{-1}{\hbar}{\rm ad}_\nu: A_\hbar \to A_\hbar\, \big\vert\, \nu \in A_\hbar \text{ is a lifting for some } x \in \big(A \bw A\big)_{\natural}  \big\}.
	\]
\end{definition}
Here, $\big(A \bw A\big)_{\natural}$ is the image of the ideal $A \bw A$ in $\HH_{0}(A)$.

\begin{definition}
	Let $(A, \ldb-,-\rdb,\bw)$ be a noncommutative Hamiltonian space.
	Let $A_{\hbar}$ be a quantization of $A$.
	The {\it noncommutative quantum gauge group of $A_\hbar$} is
	defined to be  the $A_{\hbar}$-sub-bimodule of \(\End_{\bbK[\hbar]}(A_{\hbar})\)
	generated by quantum gauge elements.
\end{definition}
Denote the noncommutative quantum gauge group by \(\mG^A_{\hbar}\).

\begin{proposition}
	Let $(A, \ldb-,-\rdb,\bw)$ be a noncommutative Hamiltonian space.
	Let $A_{\hbar}$ be a quantization of $A$.
	Then there is a canonical morphism $ \widehat{\xi}_\br: {\mG^A _{\hbar}} \to A_{\hbar} $ of $A_{\hbar}$-bimodules.
	Furthermore,
	the image of $\widehat{\xi_\br}$ coincides with  the two-sided ideal generated by
	$$
	\big\{ \ell_\hbar ([p \bw]) +\ell_{\hbar}([p\br] ) \hbar \, \vert \, [p \bw] \in (A \bw A)_\natural \big\}.
	$$
\end{proposition}
\begin{proof}
	For an arbitrary element
	$\frac{-1}{\hbar} \sum_{\al} P_\al {\rm ad}_{\ell_{\hbar}([\nu_{\alpha} \bw] )} Q_{\alpha}$ in ${\mG^A _{\hbar}}$, where $P_\al, Q_\al \in A_{\hbar}$,
	define the image  under $\widehat{\xi}_\br$ to be
	\begin{equation}\label{quant delta}
		\sum_{\al} P_\al \ell_{\hbar}([{\nu_{\alpha}} \bw ] )Q_{\alpha} +  P_\al \ell_{\hbar}([\nu_{\alpha}\br]) Q_\alpha\hbar.
	\end{equation}
	Note that $[{\nu_{\alpha}} \bw ]$ denotes the element in $\big(A \bw A\big)_{\natural}$ represented by $\nu_\al \bw \in A$.
	Then the second part of the proposition is clear.
\end{proof}

As in \cite{Zhao2021Com}, noncommutative quantum moment maps is given as follows.
\begin{definition}\label{def:  noncom quant moment map}
	Let $(A, \ldb-,-\rdb, \bw)$ be a noncommutative Hamiltonian space.
	Let $A_{\hbar}$ be a quantization of $A$.
	A {\it noncommutative quantum moment map} is defined to be a lifting of $\bw$ in $A_\hbar$.
\end{definition}

When the noncommutative quantum moment map is given, $(A_{\hbar}, \mG^A_{\hbar}, \widehat{\bw})$ is called a {\it noncommutative quantum Hamiltonian space}.
See Theorem \ref{thm: noncom quant moment map} for example.
Now,
the construction of $\widehat{\Lambda}_{A,\br} ^{\bullet}$ is clear.
It remains to prove (\ref{quantization on cohomology}).

\begin{theorem}\label{thm: quantization by a complex}
	Let $(A, \ldb-,-\rdb,\bw)$ be a noncommutative Hamiltonian space.
	Let $A_{\hbar}$ be a quantization of $A$.
	For an arbitrary $\br \in R$,
	$H^1(\widehat{\Lambda}_{A,\br} ^{\bullet})$ is a quantization of $H^1({\Lambda}_{A} ^{\bullet})$.
\end{theorem}

\begin{proof}
	By definition, to prove that 
	$$
	H^1 (\widehat{\Lambda}_{A,\br}) = \frac{A_{\hbar}}{A_{\hbar} \{\ell_\hbar ([p \bw]) +\ell_{\hbar}([p \br]) \hbar \,\vert\, [p \bw] \in (A \bw A)_\natural \} A_{\hbar} }
	$$
	is a quantization of
	\(\displaystyle{H^1(\Lambda_{A}) = \frac{A}{A \bw A}},\)
	is to prove that there exists a lifting isomorphism 
	$$
	r_\hbar: 
	Sym\big( \HH_{0}(\frac{A}{A \bw A}) \big)[\hbar]
	\longrightarrow
	\displaystyle{\frac{A_{\hbar}}{A_{\hbar} \{\ell_\hbar ([p \bw]) +\ell_{\hbar}([p \br]) \hbar\,\vert\,[p \bw] \in (A \bw A)_\natural \} A_{\hbar} }}
	$$
	of $\bbK[\hbar]$-modules such that (\ref{noncom dirac pic}) holds.

	Consider the diagram
	\begin{displaymath}
		\xymatrix{
			Sym(\HH_{0}(A))[\hbar]  \ar[d]^-{p_1}   \ar[rr]^-{\ell_\hbar}     &&   A_\hbar \ar[d]^-{p_2} \\
			Sym(\HH_{0}(B))[\hbar]   &&   B_\hbar.	
		}
	\end{displaymath}
	Here, $B$ is for the noncommutative Hamiltonian reduction $\displaystyle{H^1(\Lambda_{A}) = \frac{A}{A \bw A}}$;
	and $B_\hbar$ is for  
	\begin{equation*}
		{\frac{A_{\hbar}}{A_{\hbar} \{\ell_\hbar ([p \bw]) +\ell_{\hbar}([p \br]) \hbar\, \vert\, [p \bw] \in (A \bw A)_\natural \} A_{\hbar} }};
	\end{equation*}
	$p_1$ and $p_2$ are canonical quotient morphisms.
	By definition, an arbitrary element in $Sym(\HH_{0}(B))[\hbar]$ is of the form
	\begin{equation}\label{element1 of NQR}
		\sum_{k_1, k_2,...} c_{k_1, k_2,...} [ \bar{x}_{k_1}] \& [ \bar{x}_{k_2}] \&\cdots .
	\end{equation}
	Symbols $\bar{x}_{k_i}$ are elements in $B$
	represented by ${x}_{k_i} \in A_\hbar$,
	symbols $[ \bar{x}_{k_i}]$ are elements in $\HH_{0}(B)$,
	$c_{k_1, k_2,...}$ are polynomials with variable $\hbar$.
	Denote (\ref{element1 of NQR}) briefly by
	$\displaystyle{\sum_{\bf k }c_{\bf k} [ \bar{x}_{\bf k}]}.$
	Then the image of $\displaystyle{\sum_{\bf k }c_{\bf k} [ \bar{x}_{\bf k}]}$ under $r_\hbar$ is defined to be 
	\(\displaystyle{\sum_{k_1, k_2,...} c_{k_1, k_2,...} p_2 \circ \ell_{\hbar}\big([ {x}_{k_1}] \& [ {x}_{k_2}] \& \cdots\big)} ,\)
	and write it as
	$\sum_{\bf k }c_{\bf k}\, p_2 \circ \ell_{\hbar}\big([ {x}_{\bf k}]\big).$
	This is a well-defined  $\bbK[\hbar]$-module morphism,
	the image does not depend on representatives.
	
	$r_\hbar$ is surjective since $p_1$, $p_2$ are surjective and $\ell_{\hbar}$ is an isomorphism.
	
	Now, we show that $r_\hbar$ is injective.
	Assume  $\sum_{k_1, k_2,...} c_{k_1, k_2,...} p_2 \circ \ell_{\hbar}\big([ {x}_{k_1}] \& [ {x}_{k_2}] \& ...\big)  = 0$,
	where $\sum_{k_1, k_2,...} c_{k_1, k_2,...} [ \bar{x}_{k_1}] \& [ \bar{x}_{k_2}] \& ...  \in Sym(\HH_{0}(B))[\hbar].$
	This is equivalent to saying that
	$$
	\sum_{k_1, k_2,...} c_{k_1, k_2,...}   \ell_{\hbar}\big([ {x}_{k_1}] \& [ {x}_{k_2}] \& ...\big) \in A_{\hbar} \{\ell_\hbar ([p \bw]) +\ell_{\hbar}([p \br]) \hbar\, \vert\, [p \bw] \in (A \bw A)_\natural \}  A_{\hbar} .
	$$
	Noticed that $\ell_\hbar$ is an isomorphism and elements  of the form $[p \bw]$ are mapped to zero under $p_1$,
	then in this case, 
	$$
	\sum_{k_1, k_2,...} c_{k_1, k_2,...} [ \bar{x}_{k_1}] \& [ \bar{x}_{k_2}] \& ...  = 0.
	$$
	Consequently, $r_{\hbar}$ is an isomorphism.
	
	Since the map $r_\hbar$ is induced from $\ell_{\hbar}$, (\ref{noncom dirac pic}) holds for $r_{\hbar}$ due to $(1)$ in Proposition \ref{prop: noncom red}.
	More precisely, for any $[\bar{x}], [\bar{y}] \in \HH_{0}(B)$,
	\begin{equation*}
		\begin{split}
			r_\hbar \big(  \{ [\bar{x}], [\bar{y}] \} \big)  & = r_{\hbar}  \big( \overline{\{[x], [y]\}}\big)\\
			& = p_2 \circ \ell_{\hbar}  \big( {\{[x], [y]\}}\big)\\
			& = p_2\big( \big[\ell_{\hbar}([x]), \ell_{\hbar}([y])\big] + O(\hbar)\big)\\
			& = \big[p_2 \circ\ell_{\hbar}([x]), p_2 \circ \ell_{\hbar}([y])\big] + O(\hbar)\\
			& = \big[r_{\hbar}([\bar{x}]), r_{\hbar}([\bar{y}])\big] + O(\hbar).
		\end{split}
	\end{equation*}
	In conclusion, $H^1(\widehat{\Lambda}_{A,\br} ^{\bullet})$ is a quantization of $H^1({\Lambda}_{A} ^{\bullet})$.
\end{proof}

This theorem establishes a generalized version of the ``quantization commutes with reduction'' in noncommutative geometry, extending the observation first made in \cite{Zhao2021Com}:
\begin{equation}
	\begin{split}
		\xymatrixcolsep{4pc}
		\xymatrix{
			A_{\hbar} \ar@{-->}[r] 
			& {H^1 (\widehat{\Lambda}_{A,\br} ^\bullet)} \\
			A \ar@{~>}[u] \ar@{-->}[r]& H^1 ({\Lambda}_{A} ^\bullet).
			\ar@{~>}[u]}
	\end{split}
\end{equation}
%
%
%

\section{Quantization of preprojective algebras}\label{sec: quiver}

This section is devoted to exploring the quantization problem of quiver algebras.

\subsection{Noncommutative Hamiltonian structure on the quiver algebra}\label{subsec: necklace}

\begin{proposition}\cite[Theorem 6.3.1]{Van2008Double}
	Let $Q$ be a finite quiver.
	There is a double Poisson bracket on $\QQ$ given by the following formula:
	for any arrow $a \in Q$
	\begin{equation*}
	\ldb a, a^{\ast} \rdb = e_{s(a)} \otimes e_{t(a)},
	\ldb a^{\ast} , a \rdb = - e_{t(a)} \otimes e_{s(a)};
	\end{equation*}
	for any $f, g\in \overline{Q}  \text{ with } f \neq g^{*}$,
	\(\ldb f, g \rdb = 0.\)
\end{proposition}

The induced Lie bracket on $\HH_{0}(\QQ)$ is given as follows.
For any $a\in Q$,
$\{a,a^{*}\}=1$
and  $\{a^{*},a\}=-1$;
also,
$\{f, g\}= 0$ 
for any $f, g\in \overline{Q}  \text{ with } f \neq g^{*}$.
 By Leibniz's rule,
it is straightforward to check that
for cyclic paths $a_{1}a_{2}\cdots a_{k},b_{1}b_{2}\cdots b_{l}\in \HH_0(\QQ)$
with $a_{i},b_{j}\in \overline{Q}$, 
\begin{equation}\label{necklace bracket}
	\begin{split}
		&\{a_{1}a_{2}\cdots a_{k},b_{1}b_{2}\cdots b_{l}\}\\
		&=
		\sum_{1\leqslant i \leqslant k,\, 1\leqslant j \leqslant l}
		\{a_{i},b_{j}\}t(a_{i+1})a_{i+1}a_{i+2}\cdots a_{k}a_{1}
		\cdots a_{i-1} b_{j+1}\cdots b_{l}b_{1}\cdots b_{j-1}.
	\end{split}
\end{equation}

$\HH_{0}(\QQ) $ with the above Lie bracket is also known as 
the {\it necklace Lie algebra}.
The quiver case has been studied extensively; see \cite{CBEG2007, Van2008Double, LP1990Sem, Gin2001}.
Necessary results are summarized as the following proposition.

\begin{proposition}\cite{CBEG2007, Van2008Double}\label{prop: noncom Hamil on quiver}
	Let $Q$ be a finite quiver.
	Then the following statements hold.
	\begin{enumerate}
		\item \(\bw = \sum_{a \in Q} a a^\ast - a^{\ast} a\)
		is a noncommutative moment map for the double Poisson algebra $(\QQ, \ldb-,-\rdb)$.
		
		\item The preprojective algebra
		\(\Pi Q = \frac{\QQ}{\QQ \bw \QQ}\)
		is obtained from $\QQ$ by noncommutative Hamiltonian reduction;
		therefore, $\HH_{0}(\Pi Q)$ is  a Lie algebra, and the projection
		\(\HH_{0}(\QQ) \to \HH_{0}(\Pi Q)\)
		is a Lie algebra morphism.
	\end{enumerate}
\end{proposition}
\begin{proof}
	For the proof, see \cite[Theorem 6.3.1]{Van2008Double} and \cite[Theorem 7.2.3]{CBEG2007}.
\end{proof}

\begin{example}\label{eg: A_3}
	Consider an explicit quiver $Q$ as Figure \ref{A3}.
	Then double it, one has $\overline{Q}$ as Figure \ref{DA3}.
	\begin{figure}[H]
		\centering
		\subfigure[$Q$]{
			\includegraphics[scale=0.4]{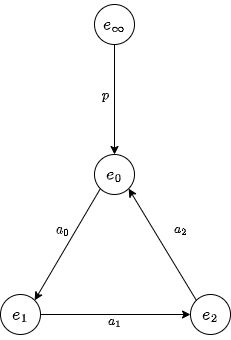}
			\label{A3}}
		\qquad\qquad\qquad
		\subfigure[$\overline Q$]{
			\includegraphics[scale=0.4]{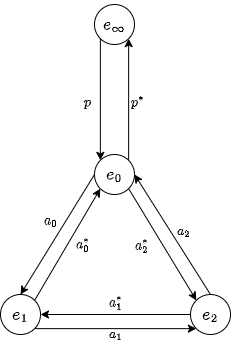}
			\label{DA3}}
		\caption{Quiver $Q$ and its doubled version $\overline{Q}$.}
	\end{figure}

	Let us calculate the necklace Lie bracket of a special cycle:
	$p^\ast a_0 ^\ast a_1^ \ast a_2 ^\ast  p$.
	\begin{equation*}
		\begin{split}
			\{p^\ast a_0 ^\ast a_1^ \ast a_2 ^\ast  p, p^\ast a_0 ^\ast a_1^ \ast a_2 ^\ast  p\}&= \{p^\ast, p\} a_0 ^\ast a_1^ \ast a_2 ^\ast  p p^\ast a_0 ^\ast a_1^ \ast a_2 ^\ast  + \{p, p^\ast\}p^\ast a_0 ^\ast a_1^ \ast a_2 ^\ast  a_0 ^\ast a_1^ \ast a_2 ^\ast p\\
			&= -a_0 ^\ast a_1^ \ast a_2 ^\ast  p p^\ast a_0 ^\ast a_1^ \ast a_2 ^\ast  + p^\ast a_0 ^\ast a_1^ \ast a_2 ^\ast  a_0 ^\ast a_1^ \ast a_2 ^\ast p\\
			&=0.
		\end{split}
	\end{equation*}
	
	Set $\Gamma = a_0 ^\ast a_1 ^ \ast a_2 ^\ast $.
	For $k, l \in \bbN$, 
	It is straightforward to check that
	\(\{p^\ast \Gamma^k p, p^\ast \Gamma^l p\} = 0.\)
	Then the trace functions $\{\Tr(p^\ast \Gamma^k p)\}_{k = 1,2,...}$ form a family of Poisson-commuting functions on $T^\ast \Rep^{Q}_{\bd}$.

	By Proposition \ref{prop: noncom Hamil on quiver}, it is clear that the noncommutative moment map is
	\begin{equation}\label{w for A_3}
		\begin{split}
			\bw & = \sum_{i=0}^{2}[a_i, a_i^\ast] + p p^\ast - p^\ast p\\
			& =  \big( - a^\ast _0 a_0 +  a_2 a_2 ^{\ast}  + p p^\ast \big)
			+ \big(a_0 a^\ast _0 - a_1^\ast a_1 \big)
			+ \big( a_1a_1^\ast -a_2 ^\ast a_2\big)
			+ \big(  - p^\ast p\big)\\
			& = \bw_0 + \bw_1 + \bw_2 + \bw_\infty.
		\end{split}
	\end{equation}
\end{example}

\subsection{Schedler's quantization}\label{subsec: schedler}

Building on \cite{Tur1991Ske}, Schedler \cite{Sch2005} constructed a PBW-deformation of the necklace Lie algebra (see also \cite{GinSch2006}).
We recall his results in this subsection.

\begin{notation}
	Let $Q$ be a finite quiver.
	Let $R$ be the semisimple ring 
	$\oplus_{i \in Q_{0}} \mathbb{K}e_{i}$.
	\begin{enumerate}
		\item Let $AH:= \overline{Q} \times \mathbb{N}$, called the space of arrows with heights. 
		
		\item Let $E_{\overline{Q},h}$ be the $\mathbb{K}$-vector space spanned by $AH$.
		
		\item Let $LH:= (T_{R}E_{\overline{Q},h})_{\natural}$, called
		the generalized cyclic path algebra with heights.
		
		\item Let $SLH[\hbar]:= Sym(LH) \otimes_{\bbK} \mathbb{K} [\hbar]$, the symmetric algebra generated by $LH$. 
	\end{enumerate}
\end{notation}

In subsequent discussions, 
all symmetric products are denoted by $\&$.
Consider the $\mathbb{K} [\hbar]$-submodule $SLH^{'}$ spanned by elements of the form
\begin{equation}
	\label{general form}
	\begin{split}
		(a_{1,1},h_{1,1})\cdots (a_{1,l_{1}},h_{1,l_{1}}) \& (a_{2,1},h_{2,1})\cdots (a_{2,l_{2}},h_{2,l_{2}})
		\\ \& \cdots  \& 
		(a_{k,1},h_{k,1})\cdots (a_{k,l_{k}},h_{k,l_{k}}) 
		\&v_{1} \& v_{2} \&\cdots \& v_{m}.
	\end{split}
\end{equation}
where the $h_{i,j}$ are all distinct, $a_{i,j}\in \overline{Q}$ and
$v_{i}\in Q_{0}.$
Let $\tilde{A}$ be the quotient of
$SLH^{'}$ where two elements in $SLH^{'}$ are identified
if and only if exchanging heights in corresponding places preserves their order.

Next, consider the $\mathbb{K}[\hbar]$-submodule $\tilde{B}$ of $\tilde{A}$ generated by the following forms.
\begin{itemize}\label{quiver skein relations}
	\item $X-X^{'}_{i,j,i^{'},j^{'}} +\hbar X^{''}_{i,j,i^{'},j^{'}}$, \\[2mm]
	where $i \neq i^{'},h_{i,j}<h_{i^{'},j^{'}},\nexists (i^{''},j^{''})\text{ with } h_{i,j} < h_{i^{''},j^{''}} < h_{i^{'},j^{'}}$;\vspace{2mm}
	\item $X-X^{'}_{i,j,i,j^{'}} + \hbar X^{''}_{i,j,i,j^{'}}$, \\[2mm]
	where $h_{i,j}<h_{i,j^{'}},\nexists (i^{''},j^{''})\text{ with } h_{i,j} < h_{i^{''},j^{''}} < h_{i^{'},j^{'}}$
\end{itemize}
In the above, $X^{'}$ and $X^{''}$ are defined as follows.
$X^{'}_{i,j,i^{'},j^{'}}$ is the same as $X$,
but with heights $h_{i,j}$ and $h_{i^{'},j^{'}}$ interchanged;
$X^{''}_{i,j,i^{'},j^{'}}$ is given by replacing the components
\[
(a_{i,1},h_{i,1})\cdots (a_{i,l_{i}},h_{i,l_{i}}) \text{ and } (a_{i^{'},1},h_{i^{'},1})\cdots (a_{i^{'},l_{i^{'}}},h_{i^{'},l_{i^{'}}})
\]
with the single component
\begin{displaymath}
	\{a_{i,j},a_{i^{'},j^{'}}\}t(a_{i,j+1})(a_{i,j+1},h_{i,j+1})\cdots (a_{i,j-1},h_{i,j-1})(a_{i^{'},j^{'}+1},h_{i^{'},j^{'}+1})\cdots (a_{i^{'},j^{'}-1},h_{i^{'},j^{'}-1}).
\end{displaymath}
$X^{'}_{i,j,i,j^{'}}$ is the same as $X$, but with heights $h_{i,j} \text{ and } h_{i,j^{'}}$ interchanged;
$X^{''}_{i,j,i,j^{'}}$ is given by replacing the component $(a_{i,1},h_{i,1})\cdots (a_{i,l_{i}},h_{i,l_{i}})$ with 
\begin{align*}
	\{(a_{i,j},a_{i,j^{'}})\} t(a_{i,j^{'}+1})(a_{i,j^{'}+1},h_{i,j^{'}+1})\cdots (a_{i,j-1},h_{i,j-1})\\
	\&
	t(a_{i,j+1})(a_{i,j+1},h_{i,j+1})\cdots (a_{i,j^{'}-1},h_{i,j^{'}-1}).
\end{align*} 

Let 
$\qneck:= {\tilde{A}}/{\tilde{B}}$.
For any $X,Y \in \qneck$, 
the product of $X$ and $Y$, denoted by $X \ast Y$, is defined to be 
``placing $Y$ above $X$". 
Throughout this work, $\qneck$ is called the {\it quantum path algebra associated with $Q$}.

\begin{proposition}\cite[Corollary 4.2]{Sch2005} \label{prop: PBW}
	Let $Q$ be a finite quiver and fix an order on the set $\{x_{i}\}$ of cyclic 
	paths in $\overline{Q}$ and idempotents in $Q_{0}$. 
	Then the projection $\tilde{A} \rightarrow Sym(\HH_{0} (\QQ))[\hbar]$ 
	obtained by forgetting the heights descends to an isomorphism 
	$$
	{\rm pr}: \qneck \rightarrow Sym(\HH_{0} (\QQ))[\hbar]
	$$
	of free $\mathbb{K}[\hbar]$-modules.
\end{proposition}

A basis of $\qneck$ as a free 
$\mathbb{K}[\hbar]$-module is given by choosing elements of the form 
(\ref{general form}) which project to the basis 
$\{[x_{i_{1}}]\&\cdots \&[x_{i_{k}}] | \text{ for any }k 
\in \mathbb{Z}_{\geq 0} \text{ and } x_{i_{1}} < x_{i_{2}} < \cdots  < x_{ i_{k}}\}$ of $Sym(\HH_{0} (\QQ))[\hbar]$.
Write
$\widehat{[x]}$ instead of $\ell_{\hbar}([x])$.
By Definition \ref{def: noncom quantization},
$\qneck$ as above is a quantization of noncommutative Poisson structure on $\QQ$.

\subsection{Noncommutative quantum reduction in the quiver setting}\label{subsec: quantum red quiver}

As an application of Theorem \ref{thm: quantization by a complex}, one has
\begin{theorem}\label{thm: noncom quant red for quiver}
	Let $Q$ be a finite quiver.
	For any $\br \in R$, 
	$$
	H^1(\widehat{\Lambda}^{\bullet}_{\QQ,  \br})  = \frac{\qneck}{\qneck   \{ \widehat{[p \bw]}+ \widehat{[p\br]}\hbar\, \vert\, [p \bw] \in (\QQ \bw \QQ)_\natural \}   \qneck}
	$$
	is a quantization of the preprojective algebra
	$
	\Pi Q.
	$
\end{theorem}

We call the algebra  $H^1(\widehat{\Lambda}^{\bullet}_{\QQ,  \br}) $
{\it quantum  preprojective algebra associated with $(Q, \br)$},
denoted by ${\Pi Q}_{\hbar, \br}$.
The above theorem implies that we have the commutative diagram
\begin{equation}
	\begin{split}
		\xymatrixcolsep{4pc}
		\xymatrix{
			\QQ_{\hbar} \ar@{-->}[r] 
			& {\Pi Q}_{\hbar, \br} \\
			\QQ \ar@{~>}[u] \ar@{-->}[r]& \Pi Q.
			\ar@{~>}[u]}
	\end{split}
\end{equation}

In  \cite{Zhao2021Com},
the quantum preprojective algebra is constructed only for $\br = 0$.
In this work, our complex formalism naturally explains the role of the parameter $\br$.
Only when the noncommutative Hamiltonian reduction is realized as a complex and the correspondence between noncommutative fields and Hamiltonians is remembered by the map $\xi$,
then correspondence between quantum fields and quantum Hamiltonians canonically contains higher-order information.
In this case,
those at order $1$ is decoded by $\br$.

\begin{example}\label{eg: quantum A_3}
	Recall that the $\overline{Q}$ is as Figure \ref{DA3}, which is 
	\begin{figure}[H]
		\centering
		{\includegraphics[scale=0.4]{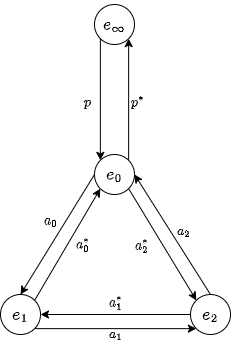}
		}
	\end{figure}
	By definition, elements in $\qneck$ can be visualized as cyclic paths with heights.
	Consider
	\begin{equation*}
		X = (a_0 ^\ast, 1)(a_1 ^\ast, 2)(a_2 ^\ast, 3)\ \text{and}\ Y = (a_2 ^\ast, 1)(a_2, 2).
	\end{equation*}
	Then
	\begin{equation*}
		\begin{split}
			[X, Y] & = (a_0 ^\ast, 1)(a_1 ^\ast, 2)(a_2 ^\ast, 3) \& (a_2 ^\ast, 4)(a_2 , 5) - (a_2 ^\ast, 1)(a_2, 2) \& (a_0 ^\ast, 3)(a_1 ^\ast, 4)(a_2 ^\ast, 5)\\
			& = (a_0 ^\ast, 2)(a_1 ^\ast, 3)(a_2 ^\ast, 4) \& (a_2 ^\ast, 1)(a_2 , 5) - (a_2 ^\ast, 1)(a_2, 2) \& (a_0 ^\ast, 3)(a_1 ^\ast, 4)(a_2 ^\ast, 5)\\
			& = - \hbar \{a_2 ^\ast, a_2\}(a_0 ^\ast, 2)(a_1 ^\ast, 3)(a_2 ^\ast, 1)\\
			& = - \hbar \{a_2 ^\ast, a_2\}(a_0 ^\ast,1)(a_1 ^\ast, 2)(a_2 ^\ast,  3 )
		\end{split}
	\end{equation*}
	On the other hand,
	\(	\{a_0 ^\ast a_1^ \ast a_2 ^\ast, a_2 ^\ast  a_2 \} = \{a_2 ^\ast, a_2\} a_0 ^\ast a_1 ^\ast a_2 ^\ast.\)
	It is clear that the projection satisfies the Dirac's picture (\ref{noncom dirac pic}) in Definition \ref{def: noncom quantization}.

	Fix elements
	\begin{equation*}
		\bw_0=  - a^\ast _0 a_0 +  a_2 a_2 ^{\ast}  + p p^\ast ,\  
		\bw_1 = a_0 a^\ast _0 - a_1^\ast a_1,\  
		\bw_2 = a_1a_1^\ast -a_2 ^\ast a_2,\ 
		\bw_\infty = - p^\ast p.
	\end{equation*}
	(\ref{w for A_3}) implies that
	noncommutative quantum gauge group is generated by
	$$
	\Big\{\frac{-1}{\hbar} {\rm ad}_{\widehat{[p \bw_{i} ]}} \, \big\vert \, p \text{ is an arbtrary cyclic path, } i = 0, 1, 2, \infty\Big\}.
	$$
\end{example}
Thus for an arbitrary parameter $\br = \sum_{i \in Q_0} r_i e_i \in R$,
noncommutative quantum reduction $H^1(\widehat{\Lambda}^{\bullet}_{\QQ,  \br})$ is the quotient algebra of $\qneck$ by the ideal generated by
\begin{equation*}
	\begin{split}
		&[- q_0a^\ast _0 a_0  +  q_0a_2 a_2 ^{\ast} + q_0 p p^\ast  ]^{\widehat{\ }} + r_0 \hbar \widehat{[q_0 ]},\ 
		[q_1 a_0 a^\ast _0 - q_1 a_1^\ast a_1 ]^{\widehat{\ }} + r_1 \hbar \widehat{[q_1 ]},\\
		&[q_2 a_1a_1^\ast  - q_2 a_2 ^\ast a_2 ]^{\widehat{\ }}  + r_2 \hbar \widehat{[q_2 ]},\ 
		[- q_\infty p^\ast p ]^{\widehat{\ }} + r_\infty \hbar \widehat{[q_\infty]}.
	\end{split}
\end{equation*}
Here, $q_i$ is an arbitrary cyclic path passing by the vertex $i$.
Since some elements are too long, we use $[-]^{\widehat{\ }}$ to represent their liftings.

\subsection{[Q, R]=0 on quiver representation spaces}\label{subsec: [Q,R]=0 on rep}

Recall that for a finite quiver $Q$ and a dimension vector $\bd$,
the quiver variety is defined to be
\(\mathcal{M}_{\bd} (Q) = \Spec \bbK[\mu^{-1}(0)]^{\GL_{\bd}(\bbK)}.\)
Here, $\mu$ is the moment map on $T^\ast \Rep^Q _\bd$:
\begin{equation*}
	\mu:  T^\ast \Rep^Q _{\bd} \to  \gl_{\bd}(\bbK)^\ast,\qquad
	 \rho \tos \tr \Big(  \big(\sum_{a \in Q}[\rho_{a}, \rho_{a^\ast}] \big) \cdot - \Big).
\end{equation*}

A quantization of a commutative Poisson algebra is defined as follows.
\begin{definition}\label{def: quant Poi}
	Suppose $S$ is a  
	commutative $\mathbb{K}$-algebra endowed 
	with Poisson bracket $\{-, -\}$. 
	A {\it  quantization of $S$} is a flat graded $\mathbb{K}[\hbar]$-algebra $S_\hbar$ $(\mathrm{deg}\hbar = 1)$
	endowed with an isomorphism $\Phi: \frac{S_{\hbar}}{\hbar S_{\hbar}}\, \rightarrow\, S$ of $\mathbb{K}$-algebras
	such that for any $a,b \in S_{\hbar}$, if we denote their images
	in $\displaystyle\frac{S_{\hbar}}{\hbar S_{\hbar}}$
	by $\overline{a},\overline{b}$, then 
	\[\Phi \Big(\overline{\frac{-1}{\hbar}(ab - ba)} \Big) 
	= \{\Phi(\overline{a}),\Phi(\overline{b})\}.
	\]
\end{definition}

It is a standard fact that the quantization of the cotangent bundle $T^\ast \Rep^Q _\bd$ is
given by the Rees algebra $\mD_{\hbar} (\Rep^Q _\bd)$ of differential operators on $\Rep^Q _\bd$.

At the quantum level, ``zero-locus defined by Hamiltonians'' is replaced by a left module algebra defined by quantum Hamiltonians (see Lu's work \cite{Lu1993} for more details).
A crucial component is the notion of quantum moment map.

\begin{definition}\label{def: quantum moment map}
	Let $G$ be an algebraic group with Lie algebra $\mathfrak{g}$.
	Let $A_{\hbar}$ be a flat graded $\mathbb{K}[\hbar]$-algebra 
	endowed with a $\mathfrak{g}$-action. The map 
	$\mu_{\hbar}:\, \mathcal{U}_{\hbar}\mathfrak{g}\, 
	\rightarrow\, A_{\hbar}$ is called a {\it quantum moment map}
	if $\mu_{\hbar}(\mathfrak{g}) \subset (A_{\hbar})_{1}$
	and for any $v \in \mathfrak{g}$, 
	$$A_{\hbar}\, \rightarrow\, A_{\hbar},\; a\mapsto\frac{-1}{\hbar}[\mu_{\hbar}(v),a]\, 
	$$ 
	is the $\mathfrak{g}$-action of 
	$v$.
\end{definition}

In the quiver case,
the quantum moment map is given by the infinitesimal action of $\GL_\bd(\bbK)$.

\begin{proposition}\label{prop: chi + habr br}
	Let $Q$ be a quiver.
	Let $\bd$ be a dimension vector.
	Then the following results hold.
	\begin{enumerate}
		\item The infinitesimal action of $\GL_\bd(\bbK)$ on $\Rep^Q _\bd$ 
		is given by the  Lie algebra homomorphism
		\(\tau:\, \gl_{\bd} (\bbK)\, \rightarrow\, 
		\mathcal{D}(\Rep^Q _\bd),\)
		which maps
		\begin{equation}\label{for: tau}
			e^{i}_{p,q}\, \mapsto\, \sum_{a\in Q,s(a)=i} \sum_{j=1}^{d_{t(a)}} [a]_{j,p} 
			\frac{\partial}{\partial (a)_{j,q}} - \sum_{a\in Q, t(a)=i}
			\sum_{j=1}^{d_{s(a)}} [a]_{q,j}\frac{\partial }{\partial (a)_{p,j}}.
		\end{equation}
		Here
		$e^{i}_{p,q}$ is the elementary matrix in the $i$-th summand of
		$\mathfrak{\gl_{\bd} (\bbK)}$.
		
		\item $\mu_{\hbar} = - \tau$
		is a quantum moment map.
		Furthermore, for any character $\chi: \gl_{\bd}(\bbK) \to \bbK$,
		$\mu_{\hbar} + \hbar \chi$ is also a quantum moment map.
	\end{enumerate}
\end{proposition}

\begin{proof}
	Since $\repqd$ is a $\GL_\bd (\bbK)$-variety, functions $\bbK[\repqd]$ is a $\GL_\bd (\bbK)$-representation.
	This implies a morphism
	\(\GL_\bd (\bbK) \to \End_{\bbK}(\bbK[\Rep^Q _\bd]);\)
	then one can define the conjugation action of $\GL_\bd (\bbK)$ on differential operators,
	i.e.
	\begin{equation*}
		g.D:=g\circ D \circ g^{-1}, \text{ for any }g \in \GL_\bd (\bbK),\ D \in \mathcal{D}(\repqd).
	\end{equation*}
	It induces $\gl_{\bd} (\bbK)$-action on the associated Rees algebra $\mathcal{D}_{\hbar}(\repqd)$:
	\begin{equation*}
		v.D:=\frac{1}{\hbar}[\tau(v), D], \text{ for any }v \in \gl_{\bd} (\bbK),\ D \in \mathcal{D}_{\hbar}(\repqd).
	\end{equation*}
	The $\tau$ is induced from Proposition \ref{prop: gauge elements}.
	Then by Definition \ref{def: quantum moment map},
	statement $(2)$ is directly induced from statement $(1)$.
\end{proof}

Here,
we adopt the following notation:
for $x\in Q$, let $[x]_{p,q}$ denote  the function $(x)_{p,q}$;
let $[x^{*}]_{p,q}$ denote the
differential operator $\displaystyle\frac{\partial}{\partial (x)_{q, p}}$. 
One can find more details in \cite[Section 3.4]{Sch2005} or \cite[Section 3.4]{Hol1999}.

Finally, Holland's result can be summarized as follows.
\begin{proposition}\cite[Proposition 2.4]{Hol1999}\label{prop: holland}
	Suppose $Q$ is a finite quiver,
	$\bd$ is a dimension vector,
	$\chi$ is a character of $\gl_{\bd} (\bbK)$ such that 
	the moment map $\mu:\, T^{*}\Rep^Q _\bd\, 
	\rightarrow \, \gl_{\bd} (\bbK)$ is a flat morphism and 
	$\chi$ vanishes on $\ker \tau$.
	Then
	$$
	\frac{ 
		\big(\mathcal{D}_{\hbar}(\Rep^Q _\bd) \big)^{\gl_{\bd} (\bbK)}}
	{ \big(\mathcal{D}_{\hbar}(\Rep^Q _\bd)(\tau - \hbar\chi)
		(\gl_{\bd} (\bbK)) \big)^{\gl_{\bd} (\bbK)}}
	$$
	is a quantization of $\mathbb{K}[\mathcal{M}_{\mathbf{d}}(Q)]$.
\end{proposition}

Therefore,
$$
\frac{ 
	\big(\mathcal{D}_{\hbar}(\Rep^Q _\bd) \big)^{\gl_{\bd} (\bbK)}}
{ \big(\mathcal{D}_{\hbar}(\Rep^Q _\bd)(\tau - \hbar\chi)
	(\gl_{\bd} (\bbK)) \big)^{\gl_{\bd} (\bbK)}}
$$
is called the {\it quantum Hamiltonian reduction of $\mD_{\hbar} (\Rep^Q _\bd)$
associated with $\chi$.}
For consistency of notations, 
let us write  
\[
\frac{ 
	\big(\mathcal{D}_{\hbar}(\Rep^Q _\bd) \big)^{\gl_{\bd} (\bbK)} }
{ \big(\mathcal{D}_{\hbar}(\Rep^Q _\bd)(\tau - \hbar\chi)
	(\gl_{\bd} (\bbK)) \big)^{\gl_{\bd} (\bbK)} }
\]
as ${\mathcal{M}_{\bd}(Q)}_{\hbar, \chi}$
and call it {\it quantum quiver variety associated with $(Q, \bd, \chi)$}.

Note that there is no general description for proper $\bd$ and $\chi$ in Proposition \ref{prop: holland},
these parameters must be determined case by case.

\begin{example}
	Let us  continue with the Example \ref{eg: A_3}.
	The dimension $\bd$ is chosen to be $d_\infty = 1$ and $d_0 = d_1 = d_2 = 2$.
	The characters of particular interest are of the form
	$\chi_{\br}=\sum_{k \in Q_{0}}\Big(\sum_{a\in Q,s(a)=k} d_{t(a)} +r_k\Big)tr_{k}$,
	see Lemma \ref{lem: q.trace preserves ideals} for details.
	To ensure that $\chi_{\br}$ vanishes on $\ker \tau$,
	one needs to solve the equation:
	$14 + 4 r_0+2 r_1 + 2r_2 = 0$.
\end{example}

At this moment, one has  two commutative diagrams:

\begin{equation*}
	\begin{split}
		\xymatrixcolsep{4pc}
		\xymatrix{
			\QQ_{\hbar} \ar@{-->}[r] 
			& {\Pi Q}_{\hbar, \br} \\
			\QQ \ar@{~>}[u] \ar@{-->}[r]& \Pi Q
			\ar@{~>}[u]}
	\end{split}
\end{equation*}
and
\begin{equation*}
	\begin{split}
		\xymatrixcolsep{4pc}
		\xymatrix{
			\mathcal{D}_{\hbar} (\Rep^Q _\bd) \ar@{-->}[r] 
			& {\mathcal{M}_{\bd}(Q)}_{\hbar, \chi} \\
			\mathbb{K}[T^{*} \Rep^Q _\bd] \ar@{~>}[u] \ar@{-->}[r]& \mathbb{K}[{\mathcal{M}_{\bd}(Q)}].
			\ar@{~>}[u]}
	\end{split}
\end{equation*}

It is natural to ask how to relate ``{\it quantization commutes with reduction}'' on quiver algebras to that on quiver representation spaces.
This will be the main goal of the rest of this section.

\subsection{Quantum trace maps}\label{subsec: q.trace}

\begin{definition}\cite[Section 3.4]{Sch2005} \label{def: q trace map}
	Suppose $Q$ is a finite quiver, $\mathbf{d}$ is a dimension vector,  
	and
	$\qneck$ is the quantum path algebra. 
	The {\it quantum trace map} $\mathrm{Tr}^q$ is a $\mathbb{K}[\hbar]$-linear map 
	from $\qneck$ to $\mathcal{D}_{\hbar}
	(\Rep^Q _\bd)$ such that 
	for any element in the form (\ref{general form}), its image is
	\begin{equation}
		d_{v1}\cdots d_{v_{m}} \sum^{d_{t(a_{i,j})}} _{\forall i,j\,k_{i,j}=1}
		\left(\prod_{h=1} ^{N} [a_{\phi ^{-1}(h)}]_{k_{\phi^{-1}(h)},
			k_{\phi^{-1}(h)+1}}\right),\label{F.1}
	\end{equation}
	where $\{ h_{i,j} \}=\{1,2,\cdots ,N\}$, 
	$\phi$ is the map such that $\phi(i,j)=h_{i,j}$.
\end{definition}

Here, $(i,j) + 1 = (i, j+1)$  with $j,\ j + 1$ taken modulo $l_i$.
It was shown in \cite[Section 3.4]{Sch2005} that $\mathrm{Tr}^q$ is independent of the choice of representatives of the elements in $\qneck$.
Furthermore, the image of a quantum trace map lies in $\gl_{\bd} (\bbK)$-invariant part.
Crucially, a quantum trace map descends to a well-defined map on noncommutative quantum reduction $H^1(\widehat{\Lambda}^{\bullet}_{\QQ, \br})$,
which is guaranteed by the following lemma.

\begin{lemma}\label{lem: q.trace preserves ideals}
	Let $Q$ be a finite quiver.
	Let $\mathbf{d}$ be a dimension vector.
	For an arbitrary $\br \in R$, there is a unique character 
	$\chi_{\br}$ of $\gl_{\bd} (\bbK)$
	such that
	\begin{align*}
		\mathrm{Tr}^{q}\Big(\qneck \{ \widehat{[p \bw]}+ \widehat{[p\br]}\hbar \, \vert\, [p \bw] &\in (\QQ \bw \QQ)_\natural \}\qneck \Big) \\
		&\subseteq \Big( \mathcal{D}_{\hbar}(\Rep^Q _\bd) (\tau - \hbar \chi_{\br}) (\gl_{\bd} (\bbK)) \Big)^{\gl_{\bd} (\bbK)}.
	\end{align*}
\end{lemma}
\begin{proof}
	By the construction of quantum path algebra $\QQ_{\hbar}$ and Definition \ref{def: q trace map},
	one only needs to check that the image of $ \{ \widehat{[p \bw]}+ \widehat{[p\br]}\hbar \, \vert \, [p \bw] \in (\QQ \bw \QQ)_\natural \}$ lies in
	\(\displaystyle{\Big( \mathcal{D}_{\hbar}(\Rep^Q _\bd) (\tau -  \hbar \chi_{\br}) 
		(\gl_{\bd} (\bbK)) \Big)^{\gl_{\bd} (\bbK)}}\)
	for some character $\chi_{\br}$.
	While the $\gl_{\bd} (\bbK)$-invariance is clear,  it remains to determine the character $\chi_{\br}$.
	
	Since $\mathrm{Tr}^{q}$ is linear, without loss of generality, 
	one can choose $\displaystyle{X=[x_{1}\cdots x_{v}\big( \sum_{a\in Q} [a, a^\ast] \big)]}
		\in (\QQ \bw \QQ)_{\natural}$ and assume $s(x_{v}) = k$. 
	Applying $\Tr^q$ to
	\begin{equation*}
		[x_{1}\cdots x_{v}\big( \sum_{a\in Q} [a, a^\ast] \big)]^{\widehat{\ }} + [x_{1}\cdots x_{v}\br ]^{\widehat{\ }}\hbar.
	\end{equation*}
	Then, one has
	\begin{align*}
		&\mathrm{Tr}^{q}\big([x_{1}\cdots x_{v}\big( \sum_{a\in Q} [a, a^\ast] \big)]^{\widehat{\ }} + [x_{1}\cdots x_{v}\br ]^{\widehat{\ }}\hbar\big)
		\\
		=&\sum_{a\in Q,t(a)=k} \sum_{l_{1},..,l_{v+2}} [x_{1}]_{l_{1},l_{2}}\cdots [x_{v}]_{l_{v},l_{v+1}} [a]_{l_{v+1},l_{v+2}} [a^{*}]_{l_{l+2},l_{1}}
		\\
		&\qquad - \sum_{a\in Q,s(a)=k} \sum_{l_{1},..,l_{v+2}} [x_{1}]_{l_{1},l_{2}}\cdots [x_{v}]_{l_{v},l_{v+1}} [a^{*}]_{l_{v+1},l_{v+2}} [a]_{l_{v+2},l_{1}}
		\\
		&\qquad +  \sum_{l_{1},..,l_{v+1}} [x_{1}]_{l_{1},l_{2}}\cdots [x_{v}]_{l_{v},l_{1}}r_k\hbar
		\\
		=&\sum_{l_{1},..,l_{v+1}}[x_{1}]_{l_{1},l_{2}}\cdots [x_{v}]_{l_{v},l_{v+1}} \sum_{a\in Q,t(a)=k} \sum_{l_{v+2}=1} ^{d_{s(a)}} [a]_{l_{v+1},l_{v+2}} [a^{*}]_{l_{v+2},l_{1}}
		\\
		&\qquad -\sum_{l_{1},..,l_{v+1}}[x_{1}]_{l_{1},l_{2}}\cdots [x_{v}]_{l_{v},l_{v+1}} \sum_{a\in Q,s(a)=k} \sum_{l_{v+2}=1}^{d_{t(a)}} \hbar \delta_{l_{1},l_{v+1}} + [a]_{l_{v+2},l_{1}} [a^*]_{l_{v+1},l_{v+2}}
		\\
		&\qquad + \sum_{l_{1},..,l_{v+1}} [x_{1}]_{l_{1},l_{2}}\cdots [x_{v}]_{l_{v},l_{v+1}}\delta_{l_1, l_{v+1}}r_k \hbar
		\\
		=&\sum_{l_{1},..,l_{v+1}}[x_{1}]_{l_{1},l_{2}}\cdots [x_{v}]_{l_{v},l_{v+1}} \Big(  \sum_{a\in Q,t(a)=k} \sum_{l_{v+2}=1} ^{d_{s(a)}}  [a]_{l_{v+1},l_{v+2}} \frac{\partial}{\partial (a)_{l_1,l_{v+2}}}
		\\
		&\qquad -\sum_{a\in Q,s(a)=k} \sum_{l_{v+2}=1} ^{d_{t(a)}} [a]_{l_{v+2},l_{1}} \frac{\partial}{\partial (a)_{l_{v+2},l_{v+1}}} - \sum_{a\in Q,s(a)=k} d_{t(a)} \delta_{l_1, l_{v+1}} \hbar + \delta_{l_1, l_{v+1}}r_k \hbar\Big)
		\\
		=&\sum_{l_{1},..,l_{v+1}}[x_{1}]_{l_{1},l_{2}}\cdots [x_{v}]_{l_{v},l_{v+1}} \Big(-\tau (e^{k}_{l_{1},l_{v+1}}) - \sum_{a\in Q, s(a)=k} (d_{t(a)} + r_k)\delta_{l_1, l_{v+1}} \hbar\Big).
	\end{align*}
	Comparing this with
	\begin{equation*}
		\tau(-e^{k}_{l_1, l_{v+1}}) - \chi_{\br }(-e^{k}_{l_1, l_{v+1}}) \hbar,
	\end{equation*}
	we obtain the character:
	\begin{equation}\label{for: r to chi}
		\chi_{\br}=\sum_{k \in Q_{0}}\Big(-\sum_{a\in Q,s(a)=k} d_{t(a)} +r_k\Big)tr_{k},
	\end{equation}
	where $tr_{k}$ denotes the trace operator on the $k$-th matrix component.
\end{proof}

Next,
we show that the noncommutative moment map fits into the Kontsevich--Rosenberg principle.
A {noncommutative quantum moment map} is a lifting of the element $\bw$ in $\qneck$.
In particular, symbol $\widehat{\bw}$ is defined to be the element:
\(\sum_{a \in Q} (a,1)(a^{*},2) - (a^{*},1)(a,2) \in \qneck.\)

\begin{theorem}\label{thm: noncom quant moment map}
	Let $Q$ be a finite quiver. Let $\bd$ be a dimension vector.
	Then the map
	\begin{equation*}
		\gl_{\bd}(\bbK) \to \mD_{\hbar}(\repqd),\ v \tos \tr([\widehat{\bw}] v)
	\end{equation*}
	is a quantum moment map for $\repqqd$.
	Furthermore, for an arbitrary $\br \in R$,
	\begin{equation*}
		\gl_{\bd}(\bbK) \to \mD_{\hbar}(\repqd),\ v \tos \tr\Big(([\widehat{\bw}] + \hbar \sum_{i \in Q_0} r_i I_i ) v \Big)
	\end{equation*}
	is also a quantum moment map.
	Here, $I_i$ is the identity matrix at $i$-th component and zero elsewhere.
\end{theorem}
\begin{proof}
	Since the trace map is linear, one only needs to prove this statement for the basis of $\gl_{\bd}(\bbK)$.
	\begin{align*}
		[\widehat{\mathbf{w}}]e^{i}_{p,q}
		&=\sum_{t(a) = i} [a][a^{*}]e^{i}_{p,q} - \sum_{s(a)=i}[a^{*}] [a] e^{i}_{p,q}\\
		&=\sum_{t(a) =i} \sum_{k,l}([a][a^{*}])_{k,l} e^{i}_{k,l} e^{i}_{p,q} 
		- \sum_{ s(a) = i} \sum_{k,l}([a^{*}] [a])_{k,l} e^{i}_{k,l} e^{i}_{p,q}\\
		&=\sum_{t(a)=i}\sum_{k} ([a][a^{*}])_{k,p} e^{i}_{k,q} - \sum_{s(a) = i} \sum_{k} ([a^{*}] [a])_{k,p} e^{i}_{k,q}.
	\end{align*}
	Then take trace on both sides, one has
	\begin{align*}
		\tr([\widehat{\bw}]e^{i}_{p,q})
		&=\tr\Big( \sum_{t(a)=i}\sum_{k} ([a][a^{*}])_{k,p} e^{i}_{k,q} - \sum_{s(a) = i} \sum_{k} ([a^{*}] [a])_{k,p} e^{i}_{k,q} \Big)\\
		&= \sum_{t(a)=i} ([a] [a^\ast])_{q,p} - \sum_{s(a) = i} ([a^\ast] [a])_{q,p} \\
		&=\sum_{t(a)=i} \sum_{l} [a]_{q, l} [a^\ast]_{l,p} - \sum_{s(a) = i} \sum_{l}[a^\ast]_{q,l} [a]_{l,p} \\
		&= \sum_{t(a)= i} \sum_{l} [a]_{q,l}\frac{\partial}{\partial (a)_{p,l}} - \sum_{s(a) = i} \sum_{l} \frac{\partial }{\partial (a)_{l,q}} [a]_{l,p}\\
		&=\sum_{t(a)= i} \sum_{l} [a]_{q,l}\frac{\partial}{\partial (a)_{p,l}} - \sum_{s(a) = i} \sum_{l} \big([a]_{l,p} \frac{\partial}{\partial (a)_{l,q}} + \hbar \delta_{p, q} \big)\\
		&=\sum_{t(a)= i} \sum_{l} [a]_{q,l}\frac{\partial}{\partial (a)_{p,l}} - \sum_{s(a) = i} \sum_{l} [a]_{l,p} \frac{\partial}{\partial (a)_{l,q}} - \hbar\sum_{s(a) = i} \sum_{l} \delta_{p, q}\\
		&=\sum_{t(a)= i} \sum_{l} [a]_{q,l}\frac{\partial}{\partial (a)_{p,l}} - \sum_{s(a) = i} \sum_{l} [a]_{l,p} \frac{\partial}{\partial (a)_{l,q}} - \hbar\sum_{s(a) = i} d_{t(a)}\delta_{p,q}.
	\end{align*}
	which is precisely 
	$(-\tau - \hbar \chi)(e^{i}_{pq})$ for some character $\chi$.
	Therefore, this proposition follows from  Proposition \ref{prop: chi + habr br}.
\end{proof}

At the end, we will explain how to relate the ``{\it quantization commutes with reduction}'' on the quiver algebra side to that on the representation space side.

\begin{proposition}\cite[Section 3.4]{Sch2005}\label{prop: [tr,quant]}
	Suppose $Q$ is a finite quiver. Then
	for any $x, y \in \HH_0(\QQ)$, one has
	$$\Phi\Big(\frac{-1}{\hbar}\big(\mathrm{Tr}^{q}(\widehat{x})\ast \Tr^q (\widehat{y})- \Tr^q (\widehat{y}) \ast \Tr^q (\widehat{x})\big)\Big)
	=\{\mathrm{Tr}(x),\mathrm{Tr}(y)\},$$
	where 
	$$\Phi: \displaystyle\frac{\mathcal{D}_{\hbar} \big( \Rep^Q _\bd \big)}{\hbar \mathcal{D}_{\hbar} \big( \Rep^Q _\bd \big)} \rightarrow \mathbb{K}[T^{*}\Rep^Q _\bd],
	[a]_{i,j}\, \mapsto\, (a)_{i,j},\ [a^{*}]_{i,j}\, \mapsto \, (a^{*})_{i,j}.$$
	In other words,  the following diagram commutes:
	\begin{equation*}\label{eq:quantizationcommwithtrace1}
		\begin{split}
			\xymatrixcolsep{4pc}
			\xymatrix{
				\QQ_{\hbar} \ar[r]^-{\mathrm{Tr}^{q}} & \mD_{\hbar}(\Rep^Q _{\bd}) \\
				\QQ \ar[r]^-{\mathrm{Tr}} \ar@{~>}[u] & \mathbb{K}[T^{*}\Rep^Q _\bd]. \ar@{~>}[u]
			}
		\end{split}
	\end{equation*}
\end{proposition}

The following example demonstrates the preservation of quantization.

\begin{example}
	Consider the $\overline{Q}$ in Example \ref{eg: A_3}
	\begin{figure}[H]
		\includegraphics[scale=0.4]{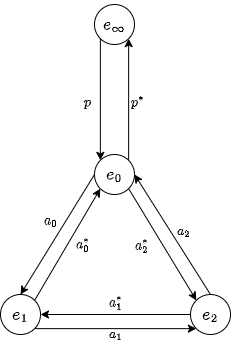}
	\end{figure}
	
	Consider 
	\begin{equation*}
		\begin{split}
			&x = a_0^\ast a_1^ \ast a_2 ^\ast\ \text{and}\ y = \ a_2^\ast  a_2 ;\\
			&\widehat{x} = (a_0 ^\ast, 1)(a_1 ^\ast, 2)(a_2 ^\ast, 3)\ \text{and}\ \widehat{y} = (a_2 ^\ast, 1)(a_2, 2).
		\end{split}
	\end{equation*}
	Then 
	\begin{equation*}
		\begin{split}
			&\frac{-1}{\hbar} \Big(\Tr^q (\widehat{x}) \ast \Tr^q (\widehat{y}) - \Tr^q (\widehat{y}) \ast \Tr^q(\widehat{x})\Big)\\
			& = \frac{-1}{\hbar}\Big(\sum_{\bf k, l} \big([a_0^\ast]_{l_1, l_2} [a_1 ^\ast]_{l_2 , l_3} [a_2 ^\ast]_{l_3, l_1} [a_2^\ast ]_{k_1 , k_2} [a_2]_{k_2, k_1} 
			- [a_2^\ast ]_{k_1 , k_2} [a_2]_{k_2, k_1}[a_0 ^\ast]_{l_1, l_2} [a_1 ^\ast]_{l_2 , l_3} [a_2 ^\ast]_{l_3, l_1}\big)\Big)\\
			& = \frac{-1}{\hbar}\Big(\sum_{\bf k, l} \big( [a_2^\ast ]_{k_1 , k_2} [a_0^\ast]_{l_1, l_2} [a_1 ^\ast]_{l_2 , l_3} [a_2 ^\ast]_{l_3, l_1}  [a_2]_{k_2, k_1} 
			- [a_2^\ast ]_{k_1 , k_2} [a_2]_{k_2, k_1}[a_0 ^\ast]_{l_1, l_2} [a_1 ^\ast]_{l_2 , l_3} [a_2 ^\ast]_{l_3, l_1}\big)\Big)  \\
			& = \frac{-1}{\hbar}\Big(\sum_{\bf k, l} \big(\hbar\delta_{k_2, l_1} \delta_{k_1, l_3}[a_2^\ast ]_{k_1 , k_2} [a_0^\ast]_{l_1, l_2} [a_1 ^\ast]_{l_2 , l_3}
			+ [a_2^ \ast ]_{k_1 , k_2} [a_2]_{k_2, k_1}   [a_0 ^\ast]_{l_1 , l_2} [a_1 ^\ast]_{l_2, l_3}  [a_2^\ast]_{l_3, l_1}\\
			&\qquad \qquad -[a_2^\ast ]_{k_1 , k_2} [a_2]_{k_2, k_1}[a_0 ^\ast]_{l_1, l_2} [a_1 ^\ast]_{l_2 , l_3} [a_2 ^\ast]_{l_3, l_1}\big)\Big)\\
			& = -\sum_{\bf k, l} \delta_{k_2, l_1} \delta_{k_1, l_3}[a_2^\ast ]_{k_1 , k_2} [a_0^\ast]_{l_1, l_2} [a_1 ^\ast]_{l_2 , l_3}\\
			& = -\sum_{\bf  l}[a_2^\ast ]_{l_3 , l_1} [a_0^\ast]_{l_1, l_2} [a_1 ^\ast]_{l_2 , l_3}.
		\end{split}
	\end{equation*}
	
	On the other hand,
	\begin{equation*}
		\begin{split}
			&\{\Tr (a_0 ^\ast a_1^ \ast a_2 ^\ast), \Tr (a_2 ^\ast  a_2 ) \}\\
			& = \{\sum_{\bf l} (a_0 ^\ast)_{l_1, l_2} (a_1 ^\ast)_{l_2, l_3} (a_2 ^\ast)_{l_3, l_1} , \sum_{\bf k} (a_2 ^\ast)_{k_1, k_2}(a_{2} )_{k_2, k_1} \}\\
			& = \sum_{\bf l, k} \{ (a_0 ^\ast)_{l_1, l_2} (a_1 ^\ast)_{l_2, l_3} (a_2 ^\ast)_{l_3, l_1},   (a_2 ^\ast)_{k_1, k_2}(a_{2} )_{k_2, k_1} \}\\
			& = \sum_{\bf l, k} \{(a_2 ^\ast)_{l_3, l_1}, (a_2)_{k_2, k_1}\} (a_0 ^\ast)_{l_1 , l_2} (a_1 ^\ast)_{l_2, l_3}  (a_2^\ast)_{k_1, k_2}\\
			& = -\sum_{\bf l}(a_0 ^\ast)_{l_1 , l_2} (a_1 ^\ast)_{l_2, l_3}  (a_2^\ast)_{l_3, l_1}.
		\end{split}
	\end{equation*}
	
	Thus, it is clear that
	\begin{equation*}
		\Phi\Big(\frac{-1}{\hbar}\big(\mathrm{Tr}^{q}(\widehat{x})\ast \Tr^q (\widehat{y})- \Tr^q (\widehat{y}) \ast \Tr^q (\widehat{x})\big)\Big)
		=\{\mathrm{Tr}(x),\mathrm{Tr}(y)\}.
	\end{equation*}
\end{example}

Analogously, one has the after-reduction case.
The argument used in the proof of \cite[Theorem 5.8]{Zhao2021Com} extends naturally to this theorem.

\begin{theorem}\label{thm: quant cond after red}
	Let $Q$ be a finite quiver.
	Let $\mathbf{d} \in \Sigma$ be a dimension vector 
	such that the moment map $\mu$ is flat. Then for any $\br \in  R$ and any $x, y \in \HH_{0}(\Pi Q)$, 
	one has
	$$\Phi \Big( \frac{-1}{\hbar}\big(\mathrm{Tr}^{q}(\widehat{x}) \ast \Tr^q (\widehat{y}) - \Tr^q (\widehat{y}) \ast \Tr^q (\widehat{x})\big) \Big)
	= \{\mathrm{Tr}(x),\mathrm{Tr}(y)\}.
	$$ 
	In other words, 
	the following diagram is commutative:
	\begin{equation*}\label{eq:quantizationcommwithtrace2}
		\begin{split}
			\xymatrixcolsep{4pc}
			\xymatrix{
				\Pi Q_{\hbar, \br} 
				\ar[r]^-{\mathrm{Tr}^{q}} & \mathcal{M}_{\bd}(Q)_{\hbar, \chi_{\br}} \\
				\Pi Q \ar[r]^-{\mathrm{Tr}} \ar@{~>}[u]& 
				\mathbb{K}[\mathcal{M}_{\mathbf{d}}(Q)]. \ar@{~>}[u]
			}
		\end{split}
	\end{equation*} 
\end{theorem}

In conclusion, one has the following commutative cubic:
\begin{equation*}
	\begin{split}
		\xymatrixrowsep{0.8pc}
		\xymatrixcolsep{1.2pc}
		\xymatrix{
			\qneck \ar@{-->}[rd] \ar[rr]^-{\mathrm{Tr}^{q}}  && 
			\mD_{\hbar}(\Rep^Q _{\bd}) \ar@{-->}[rd] \\
			&\Pi Q_{\hbar, \br} \ar[rr]^{\mathrm{Tr}^{q}}  
			&& 
			\mathcal{M}_{\bd}(Q)_{\hbar, \chi_{\br}} \\
			\mathbb{K}\overline{Q} \ar@{~>}[uu] \ar@{-->}[rd] \ar'[r][rr]^-{\mathrm{Tr}}&& 
			\mathbb{K}[T^\ast \Rep^Q _\bd] \ar@{~>}'[u][uu] \ar@{-->}[rd] \\
			& \Pi Q \ar@{~>}[uu] \ar[rr]^{\mathrm{Tr}} 
			&& \mathbb{K}[\mathcal{M} _{\mathbf{d}}(Q)] \ar@{~>}[uu]
		}
	\end{split}
\end{equation*}

\subsection{The deformed case}\label{subsec: deformed case}

Let $Q$ be a finite quiver. Let $(\QQ, \ldb-,-\rdb)$ be the double Poisson algebra given in Section \ref{subsec: necklace}.
For an arbitrary $\lambda = \sum_{i \in Q_0} \lambda_i e_i \in R = \oplus_{i \in Q_{0}} \bbK e_i$, one can easily check that
\(\bw - \lambda = \sum_{a \in Q} [a , a^\ast] - \lambda\)
is a noncommutative moment map for $(\QQ, \ldb-,-\rdb)$.
In other words,
the deformed preprojective algebra
\(\Pi^\lambda Q = \frac{\QQ}{\QQ (\bw - \lambda) \QQ}\)
is a noncommutative Hamiltonian reduction of $\QQ$.

Then,
as a corollary of Theorem \ref{thm: quantization by a complex},  one has
\begin{corollary}
	Let $Q$ be a finite quiver. For arbitrary $\lambda,\ \br \in R$, one has a quantization of the deformed preprojective algebra
	\begin{equation*}
		\Pi^\lambda Q_{\hbar, \br} = \frac{\qneck}{\qneck   \{ \widehat{[p \bw]}-\widehat{[p\lambda]}+ \widehat{[p\br]}\hbar\, \vert\, [p \bw] \in (\QQ \bw \QQ)_\natural \}   \qneck}
	\end{equation*}
\end{corollary} 

The quantum trace map descends to a well-defined map on the quantum deformed preprojective algebra $\Pi^\lambda Q_{\hbar, \br}$.
\begin{corollary}
	Let $Q$ be a finite quiver.
	Let $\mathbf{d}$ be a dimension vector.
	For  any $\lambda,\ \br \in R$
	with $\sum_{i \in Q_0} \lambda_i d_i = 0$, there is a unique character 
	$\chi_{\br}$ of $\gl_{\bd} (\bbK)$
	such that
	\begin{equation*}
		\mathrm{Tr}^{q}\Big(\qneck   \{ \widehat{[p \bw]}-\widehat{[p\lambda]}+ \widehat{[p\br]}\hbar\, \vert\, [p \bw] \in (\QQ \bw \QQ)_\natural \}   \qneck \Big)
	\end{equation*}
	is contained in
	\begin{equation*}
		\Big( \mathcal{D}_{\hbar}(\Rep^Q _\bd) (\tau + \sum_{i \in Q_0} \lambda_i tr_i - \hbar \chi_{\br}) 
		(\gl_{\bd} (\bbK)) \Big)^{\gl_{\bd} (\bbK)}.
	\end{equation*}
\end{corollary}
\begin{proof}
	The proof is the same as Lemma \ref{lem: q.trace preserves ideals} by noticing that
	\begin{align*}
		&\mathrm{Tr}^{q}\big( [ x_{1}\cdots x_{v}\big( \sum_{a\in Q} [a, a^\ast] \big) ]^{\widehat{\ }} - [x_{1}\cdots x_{v} \lambda]^{\widehat{\ }} + [x_{1}\cdots x_{v}\br ]^{\widehat{\ }} \hbar \big)\\
		=&\sum_{l_{1},..,l_{v+1}}[x_{1}]_{l_{1},l_{2}}\cdots [x_{v}]_{l_{v},l_{v+1}} \Big(-\tau (e^{k}_{l_{1},l_{v+1}}) -\lambda_k \delta_{l_1, l_{v+1}}- \sum_{a\in Q, s(a)=k} (d_{t(a)} + r_k)\delta_{l_1, l_{v+1}} \hbar\Big).
	\end{align*}
	Then, compare it with 
	\begin{equation*}
		\tau(-e^{k}_{l_1, l_{v+1}}) + \sum_{i \in Q_0} \lambda_i tr_i(-e^{k}_{l_1, l_{v+1}}) - \chi_{\br }(-e^{k}_{l_1, l_{v+1}}) \hbar,
	\end{equation*}
	the corollary holds.
\end{proof}

According to \cite[Proposition 2.4]{Hol1999},
Proposition \ref{prop: holland} holds in the deformed preprojective algebra case;
then follow analysis in Section \ref{subsec: q.trace},
one also has the commutative cubic for deformed preprojective algebras
\begin{equation*}
	\begin{split}
		\xymatrixrowsep{0.8pc}
		\xymatrixcolsep{1.2pc}
		\xymatrix{
			\qneck \ar@{-->}[rd] \ar[rr]^-{\mathrm{Tr}^{q}}  && 
			\mD_{\hbar}(\Rep^Q _{\bd}) \ar@{-->}[rd] \\
			&\Pi^\lambda Q_{\hbar, \br} \ar[rr]^{\mathrm{Tr}^{q}}  
			&& 
			\mathcal{M}^\lambda _{\bd}(Q)_{\hbar, \chi_{\br}} \\
			\mathbb{K}\overline{Q} \ar@{~>}[uu] \ar@{-->}[rd] \ar'[r][rr]^-{\mathrm{Tr}}&& 
			\mathbb{K}[T^\ast \Rep^Q _\bd] \ar@{~>}'[u][uu] \ar@{-->}[rd] \\
			& \Pi^\lambda Q \ar@{~>}[uu] \ar[rr]^{\mathrm{Tr}} 
			&& \mathbb{K}[\mathcal{M}^\lambda _{\mathbf{d}}(Q)] \ar@{~>}[uu]
		}
	\end{split}
\end{equation*}
for proper choice of $\lambda,\ \br,\bd $
such that assumptions in Proposition \ref{prop: holland} hold in the deformed setting.
Here,
$\mathbb{K}[\mathcal{M}^\lambda _{\mathbf{d}}(Q)] = \bbK[\mu^{-1} (\lambda)]^{\GL_{\bd} (\bbK)}.$

\bibliographystyle{plain}
\bibliography{bib_noncomquantred}

\begin{thebibliography}{10}

\bibitem{BL2021Eti}
R.~Bezrukavnikov and I.~Losev.
\newblock Etingof's conjecture for quantized quiver varieties.
\newblock {\em Invent. Math.}, 223(3):1097--1226, 2021.

\bibitem{BLPB2016II}
T.~Braden, A.~Licata, N.~Proudfoot, and B.~Webster.
\newblock Quantizations of conical symplectic resolutions {II}: category
  $\mathcal{O}$ and symplectic duality.
\newblock {\em Ast\'{e}risque}, (384):75--179, 2016.
\newblock with an appendix by I. Losev.

\bibitem{BPW2016Qua}
T.~Braden, N.~Proudfoot, and B.~Webster.
\newblock Quantizations of conical symplectic resolutions {I}: local and global
  structure.
\newblock {\em Ast\'{e}risque}, (384):1--73, 2016.

\bibitem{CG2010Rep}
Neil Chriss and Victor Ginzburg.
\newblock {\em Representation theory and complex geometry}.
\newblock Modern Birkh\"{a}user Classics. Birkh\"{a}user Boston, Ltd., Boston,
  MA, 2010.
\newblock Reprint of the 1997 edition.

\bibitem{CBEG2007}
W.~Crawley-Boevey, P.~Etingof, and V.~Ginzburg.
\newblock Noncommutative geometry and quiver algebras.
\newblock {\em Adv. Math.}, 209(1):274--336, 2007.

\bibitem{EGGO2007HC}
P.~Etingof, W.~L. Gan, V.~Ginzburg, and A.~Oblomkov.
\newblock Harish-{C}handra homomorphisms and symplectic reflection algebras for
  wreath-products.
\newblock {\em Publ. Math. Inst. Hautes \'{E}tudes Sci.}, (105):91--155, 2007.

\bibitem{EG2002}
P.~Etingof and V.~Ginzburg.
\newblock Symplectic reflection algebras, {C}alogero-{M}oser space, and
  deformed {H}arish-{C}handra homomorphism.
\newblock {\em Invent. Math.}, 147(2):243--348, 2002.

\bibitem{Gin2001}
V.~Ginzburg.
\newblock Non-commutative symplectic geometry, quiver varieties, and operads.
\newblock {\em Math. Res. Lett.}, 8(3):377--400, 2001.

\bibitem{Gin2005non}
V.~Ginzburg.
\newblock Lectures on noncommutative geometry.
\newblock {\em arXiv: Algebraic Geometry}, 2005.

\bibitem{GinSch2006}
V.~Ginzburg and T.~Schedler.
\newblock Moyal quantization and stable homology of necklace {L}ie algebras.
\newblock {\em Mosc. Math. J.}, 6(3):431--459, 587, 2006.

\bibitem{Hol1999}
M.~P. Holland.
\newblock Quantization of the {M}arsden-{W}einstein reduction for extended
  {D}ynkin quivers.
\newblock {\em Ann. Sci. \'{E}cole Norm. Sup. (4)}, 32(6):813--834, 1999.

\bibitem{Kontsevich}
M.~Kontsevich.
\newblock Formal (non)commutative symplectic geometry.
\newblock In {\em The {G}elfand {M}athematical {S}eminars, 1990--1992}, pages
  173--187. Birkh\"{a}user Boston, Boston, MA, 1993.

\bibitem{KonRos2000}
M.~Kontsevich and A.~L. Rosenberg.
\newblock Noncommutative smooth spaces.
\newblock In {\em The {G}elfand {M}athematical {S}eminars, 1996--1999}, Gelfand
  Math. Sem., pages 85--108. Birkh\"{a}user Boston, Boston, MA, 2000.

\bibitem{LP1990Sem}
L.~Le~Bruyn and C.~Procesi.
\newblock Semisimple representations of quivers.
\newblock {\em Trans. Amer. Math. Soc.}, 317(2):585--598, 1990.

\bibitem{Los2021Loc}
I.~Losev.
\newblock Localization theorems for quantized symplectic resolutions.
\newblock 2021.

\bibitem{Lu1993}
J.-H. Lu.
\newblock Moment maps at the quantum level.
\newblock {\em Comm. Math. Phys.}, 157(2):389--404, 1993.

\bibitem{Sch2005}
T.~Schedler.
\newblock A {H}opf algebra quantizing a necklace {L}ie algebra canonically
  associated to a quiver.
\newblock {\em Int. Math. Res. Not.}, (12):725--760, 2005.

\bibitem{Soi2019Lec}
Alexander Soibelman.
\newblock Lecture notes on quiver representations and moduli problems in
  algebraic geometry.
\newblock {\em arXiv: Algebraic Geometry}, 2019.

\bibitem{Tur1991Ske}
V.~G. Turaev.
\newblock Skein quantization of {P}oisson algebras of loops on surfaces.
\newblock {\em Ann. Sci. \'{E}cole Norm. Sup. (4)}, 24(6):635--704, 1991.

\bibitem{Van2008Double}
M.~Van~den Bergh.
\newblock Double {P}oisson algebras.
\newblock {\em Trans. Amer. Math. Soc.}, 360(11):5711--5769, 2008.

\bibitem{Zhao2021Com}
H.~Zhao.
\newblock Commutativity of quantization and reduction for quiver
  representations.
\newblock {\em Math. Z.}, 301(4):3525--3554, 2022.

\end{thebibliography}

\end{document}